# Continued Fractions and Probability Estimations in Shor's Algorithm: A Detailed and Self-Contained Treatise

Johanna Barzen[0000-0001-8397-7973] and Frank Leymann[0000-0002-9123-259X]

University of Stuttgart, IAAS, Universitätsstr. 38, 70569 Stuttgart, Germany
{firstname.lastname}@iaas.uni-stuttgart.de

**Abstract.** The algorithm of Shor for prime factorization is a hybrid algorithm consisting of a quantum part and a classical part. The main focus of the classical part is a continued fraction analysis. The presentation of this is often short, pointing to text books on number theory. In this contribution, we present the relevant results and proofs from the theory of continued fractions in detail (even in more detail than in text books) filling the gap to allow a complete comprehension of Shor's algorithm. Similarly, we provide a detailed computation of the estimation of the probability that convergents will provide the period required for determining a prime factor.

## 1. Introduction

The algorithm of Shor [7] for prime factorization is generally considered as a major milestone and a breakthrough in quantum computing: it solves a practically very relevant problem (which is, e.g., an underpinning of cryptography) with an exponential speedup compared to classical methods.

The overall algorithm is hybrid, consisting of classical computations and a quantum computation. The classical computations are computing greatest common divisors with the Euclidian algorithm, and perform a continuous fraction analysis. A detailed discussion of the latter is one of the two foci of this contribution (see part I).

The quantum part mainly consists of (i) creating an entangled state based on an oracle $f$ computing a modular exponentiation, (ii) performing a quantum Fourier transform (QFT) on this state, and (iii) measuring it. The oracle produces the following state:

$$|a\rangle\,|b\rangle = \frac{1}{\sqrt{N}} \sum_{x=0}^{N-1} |x\rangle\,|f(x)\rangle \tag{1}$$

After applying the quantum Fourier transform and a measurement, the first part (i.e. the $|a\rangle$-part) of the quantum register is in state

$$|y\rangle = \frac{1}{\sqrt{NA}} \sum_{j=0}^{A-1} \omega_N^{jpy} |y\rangle \tag{2}$$

The measured value y can then be used with high probability (see section 10, theorem 60) to compute the period of the modular exponentiation function $f$ by analyzing the convergents of a continued fraction (see section 10.1) and finally, based on the period, a prime factor (see section 8.1).

## 1.1. Structure of the Article

The article is structured as follows: in part I we cover all details about continued fractions that are required to comprehend the corresponding aspect of Shor's algorithm.

Section 2 defines the notion of a continued fraction, gives examples of how to compute the continued fraction representation of a rational number, and how to compute the number that a continued fraction (and thus convergents) represent.

Convergents as the fundamental tool in the theory of continued fractions are detailed in section 3: after defining the term, basic theorems about convergents like the recursion theorem, two sign theorems, monotony properties, convergent comparison, nesting of a number by its convergents, and several distance estimations are proven.

Next, the brief section 4 presents infinite regular continued fractions to represent non-rational numbers. A corresponding algorithms is provided to compute such continued fractions.

Section 5 gives several upper bounds and lower bounds for the difference between a number and its convergents. Exploiting one of these bounds, the convergence of the convergents of an infinite regular continued fraction of a number to this number is proven. Semiconvergents are defined and corresponding monotony properties are given.

Best approximations of a real number are introduced in section 6. It is proven, that best approximations of the second kind are convergents and vice versa (theorem of Lagrange). Best approximations of the first kind are proven to be convergents or semiconvergents (another theorem of Lagrange). Finally, the theorem of Legendre is presented which is the main result about continued fractions required by the Shor algorithm: it allows to imply that a given fraction is a convergent of another number.

Part II is devoted to estimating the probability that convergents can be used to compute periods, i.e. that Legendre's theorem can be applied.

At the beginning of part II, section 7 proves a lower bound and an upper bound for the secant lengths of the unit circle. This estimation is central for estimating the before mentioned probability.

Section 8 contains many different estimations of parameters that appear in the measurement result of the Shor algorithm. In 8.1 we remind the very basics of modular arithmetics, relate this to group theory, and use the Lagrange theorem from group theory to prove that the period of the modular exponentiation function in Shor's algorithm is less than the number to be factorized (lemma 41). Intervals of consecutive multiples of the period are studied in section 8.2: it is shown that multiples of N are sparsely scattered across these intervals (note 49). This implies that measurement results are somehow centered around multiples of $N/p$ (corollary 52). The cardinality of arguments in the superposition that build the pre-image of a certain $f(x)$ is estimated in section 8.3. Section 8.4 proves bounds of phases of amplitudes relevant for computing the probability of measurement results as a geometric sum.

Finally, section 9 computes this probability: It is proven that a measurement result is close to a multiple of $N/p$ with probability of approximately $4/\pi^2$ (lemma 57).

Section 10 shows that this measurement result fulfills the assumption of the Legendre theorem (theorem 59). Thus, by computing convergents the period can be determined (theorem 60 and section 10.1).

A brief conclusion and discussion of related work ends this contribution with section 11.

# Part I: Continued Fractions

## 2. Definition of Continued Fractions and Their Computation

We define the notion of continued fractions and give an example of how to compute them.

**Definition 1**: An expression of the form

$$a_0 + \cfrac{b_1}{a_1 + \cfrac{b_2}{a_2 + \cfrac{b_3}{\ddots}}} \tag{3}$$

with $a_i, b_i \in \mathbb{C}$ is called an *infinite* continued fraction.

If in this expression, it is $b_i = 1$ for all i, $a_0 \in \mathbb{Z}$, and $a_i \in \mathbb{N}$ for i≥1, the expression is called *regular* continued fraction.

A *finite* regular continued fraction (simply called *continued fraction*) satisfies in addition the condition $\exists\, N \in \mathbb{N}\ \forall\, k \in \mathbb{N} : a_{N+k} = 0$ (convention: "1/0 = 0") .

A continued fraction is, thus, the following expression:

$$[a_0; a_1, \cdots, a_N] \overset{\text{def}}{=} a_0 + \cfrac{1}{a_1 + \cfrac{1}{\ddots + \cfrac{1}{a_{N-1} + \cfrac{1}{a_N}}}} \tag{4}$$

□

A continued fraction of a rational number a/b is computed as follows: the integer part $\lfloor a/b \rfloor$ becomes $a_0 \in \mathbb{Z}$ leaving the non-negative rational remainder $x_1/y_1 \in \mathbb{Q}$. The latter is now written as $1/(y_1/x_1)$ resulting in

$$a_0 + \cfrac{1}{\left(\cfrac{y_1}{x_1}\right)}$$

Next, the integer part $\lfloor y_1/x_1 \rfloor$ becomes $a_1$, leaving a rational remainder that is treated as before. This processing stops until the rational remainder is zero. Figure 1 gives an example of the processing.

$$\frac{67}{47} = 1 + \frac{20}{47} = 1 + \cfrac{1}{\cfrac{47}{20}} = 1 + \cfrac{1}{2 + \cfrac{7}{20}} = 1 + \cfrac{1}{2 + \cfrac{1}{\cfrac{20}{7}}}$$

$$= 1 + \cfrac{1}{2 + \cfrac{1}{2 + \cfrac{6}{7}}} = 1 + \cfrac{1}{2 + \cfrac{1}{2 + \cfrac{1}{\cfrac{7}{6}}}} = 1 + \cfrac{1}{2 + \cfrac{1}{2 + \cfrac{1}{1 + \cfrac{1}{6}}}}$$

$$\Rightarrow \frac{67}{47} = [1;2,2,1,6]$$

**Fig. 1**. Example of a straightforward computation of a continued fraction.

Beside this straightforward proceeding to compute continues fractions, the well-known Euclidian algorithm can be used for this purpose too. Figure 2 gives a corresponding example; it should be self-descriptive.

$$43 = \mathbf{2} \times 19 + 5$$
$$19 = \mathbf{3} \times 5 + 4$$
$$5 = \mathbf{1} \times 4 + 1$$
$$4 = \mathbf{4} \times 1 + 0$$

$$\Rightarrow \frac{43}{19} = [2;3,1,4]$$

**Fig. 2**. Using the Euclidian algorithm to compute a continued fraction.

Formally, a continued fraction can always be reduced such that its last element is greater than or equal to 2.

**Note 2**

Let $[a_0; a_1, \ldots, a_N]$ be a continued fraction. Then:

$$[a_0; a_1, \ldots, a_N] = \left[a_0; a_1, \ldots, a_{N-1} + \frac{1}{a_N}\right] \tag{5}$$

Especially, it can always be achieved that a continued fraction $[a_0; a_1, \ldots, a_N]$ satisfies $a_N \geq 2$.

Proof: The following simple computation proves the first claim:

$$\begin{aligned}
[a_0; a_1, \ldots, a_N] &= a_0 + \cfrac{1}{a_1 + \cfrac{1}{\ddots + \cfrac{1}{a_{N-1} + \cfrac{1}{a_N}}}} \\
&= a_0 + \cfrac{1}{a_1 + \cfrac{1}{\ddots + \cfrac{1}{\left(a_{N-1} + \cfrac{1}{a_N}\right)}}} \\
&= \left[a_0; a_1, \ldots, a_{N-1} + \frac{1}{a_N}\right]
\end{aligned}$$

Furthermore, if $a_N = 1$ in $[a_0; a_1, \ldots, a_N]$ then $a_{N-1} + 1/a_N \geq 2$. This is because by definition $a_k \geq 1$ for $1 \leq k \leq N$.

■

Equation (5) implies a straightforward way to compute the value represented by a continued fraction $[a_0; a_1, \ldots, a_N]$: see figure 3.

$$\begin{aligned}
[2; 3, 1, 4] &= [2; 3, 1 + 1/4] = [2; 3, 5/4] = [2; 3 + 1/(5/4)] \\
&= [2; 3 + 4/5] = [2; 19/5] = [2 + 1/(19/5)] = [2 + 5/19] \\
&= [43/19] = \frac{43}{19}
\end{aligned}$$

**Fig. 3**. Computing the value of a continued fraction based on equation (5).

## 3. Convergents

Next, we define the "workhorses" of the theory of continued fractions.

**Definition 3**: $[a_0; a_1, \ldots, a_m]$ is called *m-th convergent* of the continued fraction $[a_0; a_1, \ldots, a_N]$ for 0 ≤ m ≤ N, or m-th convergent of the infinite regular continued fraction $[a_0; a_1, \ldots]$. □

Convergents can be computed recursively based on the following theorem.

**Theorem 4** (*Recursion Theorem*)
Define

- $p_0 = a_0$,
- $p_1 = a_1 a_0 + 1$ and
- $p_n = a_n p_{n-1} + p_{n-2}$ for n ≥ 2,

and define

- $q_0 = 1$,
- $q_1 = a_1$ and
- $q_n = a_n q_{n-1} + q_{n-2}$ for n ≥ 2.

Then, for every convergent $[a_0; a_1, \ldots, a_n]$ it is :

$$[a_0; a_1, \ldots, a_n] = \frac{p_n}{q_n} \tag{6}$$

<u>Proof</u> (by induction): Let n= 0, 1: Then $[a_0] = a_0 = \frac{p_0}{q_0}$ and $[a_0; a_1] = a_0 + \frac{1}{a_1} = \frac{a_0 a_1 + 1}{a_1} = \frac{p_1}{q_1}$.

Induction hypothesis: $[a_0; a_1, \ldots, a_n] = \frac{p_n}{q_n} = \frac{a_n p_{n-1} + p_{n-2}}{a_n q_{n-1} + q_{n-2}}$.

Induction step n ⟶ n+1: According to note 2, it is

$$[a_0; a_1, \ldots, a_n, a_{n+1}] = \left[a_0; a_1, \ldots, a_n + \frac{1}{a_{n+1}}\right],$$

and the last continued fraction has n elements, i.e. the induction hypothesis applies:

$$\left[a_0; a_1, \ldots, a_n + \frac{1}{a_{n+1}}\right] = \frac{\left(a_n + \frac{1}{a_{n+1}}\right) p_{n-1} + p_{n-2}}{\left(a_n + \frac{1}{a_{n+1}}\right) q_{n-1} + q_{n-2}} = \frac{\frac{a_n a_{n+1} + 1}{a_{n+1}} p_{n-1} + p_{n-2}}{\frac{a_n a_{n+1} + 1}{a_{n+1}} q_{n-1} + q_{n-2}}$$

$$= \frac{\left(a_n a_{n+1} + 1\right) p_{n-1} + a_{n+1} p_{n-2}}{\left(a_n a_{n+1} + 1\right) q_{n-1} + a_{n+1} q_{n-2}}$$

$$= \frac{a_{n+1}\left(a_n p_{n-1} + p_{n-2}\right) + p_{n-1}}{a_{n+1}\left(a_n q_{n-1} + q_{n-2}\right) + q_{n-1}} \overset{(A)}{=} \frac{a_{n+1} p_n + p_{n-1}}{a_{n+1} q_n + q_{n-1}}$$

$$\overset{(B)}{=} \frac{p_{n+1}}{q_{n+1}}$$

Here, (A) is valid because of the induction hypothesis, and (B) is the definition of $p_{n+1}$ and $q_{n+1}$.

■

The recursion theorem implies the often used

**Corollary 5**
Numerators and denominators of convergents of a continued fraction $\left[a_0; a_1, \ldots, a_N\right]$ with $a_0 \geq 0$ are strictly monotonically increasing:

$$p_n > p_{n-1} \text{ and } q_n > q_{n-1} \text{ for all } n \in \mathbb{N}.$$

Proof (by induction):
Let n=1: By definition, $p_0 = a_0, p_1 = a_1 a_0 + 1$. Because $a_i \geq 1$ for i ≥ 1, and $a_0 \geq 0$, it is $p_1 > p_0 \geq 0$. Similarly, $q_1 > q_0 > 0$

Now, $p_n = a_n p_{n-1} + p_{n-2}$ and $q_n = a_n q_{n-1} + q_{n-2}$ for n ≥ 2. With $a_n \geq 1$ by definition, and $p_{n-1} > p_{n-2}$ (≥1) as well as $q_{n-1} > q_{n-2}$ (≥1) by induction hypothesis, the claim follows.

■

The next theorem is about the sign of a combination of the numerators and denominators of consecutive convergents of a continued fraction.

**Theorem 6** (*Sign Theorem*)
For $\left[a_0; a_1, \ldots, a_n\right] = \frac{p_n}{q_n}$ the following holds:

$$p_n q_{n-1} - p_{n-1} q_n = (-1)^{n-1} \tag{7}$$

Proof (by induction): For n = 1 it is
$p_1 q_0 - p_0 q_1 = (a_1 a_0 + 1) \cdot 1 - a_0 \cdot a_1 = 1 = (-1)^0$.
Induction step n ⟶ n+1:

$$\begin{aligned} p_{n+1}q_n - p_nq_{n+1} &= \left(a_{n+1}p_n + p_{n-1}\right) q_n - p_n \left(a_{n+1}q_n + q_{n-1}\right) \\ &= a_{n+1}p_nq_n + p_{n-1}q_n - p_na_{n+1}q_n - p_nq_{n-1} \\ &= p_{n-1}q_n - p_nq_{n-1} = -\left(p_nq_{n-1} - p_{n-1}q_n\right) \\ &\overset{(A)}{=} -(-1)^{n-1} = (-1)^n \end{aligned}$$

(A) uses the induction hypothesis.

■

In case the numerators and denominators stem from the n-th convergent and the (n−2)-nd convergent, the last n-th element of the convergent gets part of the equation.

**Theorem 7** (*Second Sign Theorem*)
For $\left[a_0; a_1, \ldots, a_n\right] = \frac{p_n}{q_n}$ the following holds:

$$p_nq_{n-2} - p_{n-2}q_n = (-1)^n a_n \tag{8}$$

Proof: It is $p_n = a_np_{n-1} + p_{n-2}$ and $q_n = a_nq_{n-1} + q_{n-2}$.
Multiplying the first equations by $q_{n-2}$ and the second equation by $p_{n-2}$ results in $q_{n-2}p_n = q_{n-2}a_np_{n-1} + q_{n-2}p_{n-2}$ and $p_{n-2}q_n = p_{n-2}a_nq_{n-1} + p_{n-2}q_{n-2}$. Next, both equations are subtracted:

$$\begin{aligned} p_nq_{n-2} - p_{n-2}q_n &= q_{n-2}a_np_{n-1} + q_{n-2}p_{n-2} - p_{n-2}a_nq_{n-1} - p_{n-2}q_{n-2} \\ &= q_{n-2}a_np_{n-1} - p_{n-2}a_nq_{n-1} \\ &= a_n\left(p_{n-1}q_{n-2} - p_{n-2}q_{n-1}\right) \\ &\overset{(A)}{=} (-1)^n a_n \end{aligned}$$

where (A) is implied by the sign theorem (6) and considering $(-1)^{n-2} = (-1)^n$.

■

The sign theorem yields immediately the important

**Corollary 8**
Numerator and denominator of a convergent are co-prime.

Proof: Let t be a divisor of $p_n$ and $q_n$, i.e. $t \mid p_n$ and $t \mid q_n$. Then $t \mid (p_nq_{n-1} - p_{n-1}q_n)$, but $(p_nq_{n-1} - p_{n-1}q_n) = (-1)^{n-1}$ according to the sign theorem. Thus, t = ±1.

■

Convergents can be represented as a sum of fractions with alternating sign and whose denominators consists of products of two consecutive denominators from the recursion theorem.

**Theorem 9** (*Representation as a Sum*)
Each convergent can be represented as a sum:

$$[a_0; a_1, \dots, a_n] = a_0 + \frac{1}{q_1 q_0} - \frac{1}{q_2 q_1} + \cdots + (-1)^{n-1} \frac{1}{q_n q_{n-1}} \tag{9}$$

Proof: Let $[a_0; a_1, \dots, a_n] = \frac{p_n}{q_n}$. Since $-\frac{p_i}{q_i} + \frac{p_i}{q_i} = 0$, we can write

$$[a_0; a_1, \dots, a_n] = \frac{p_n}{q_n} - \frac{p_{n-1}}{q_{n-1}} + \frac{p_{n-1}}{q_{n-1}} - \frac{p_{n-2}}{q_{n-2}} + \frac{p_{n-2}}{q_{n-2}} - \cdots + \frac{p_1}{q_1} - \frac{p_0}{q_0} + \frac{p_0}{q_0}$$

Computing the differences results in

$$[a_0; a_1, \dots, a_n] = \frac{p_n q_{n-1} - q_n p_{n-1}}{q_n q_{n-1}} + \frac{p_{n-1} q_{n-2} - q_{n-1} p_{n-2}}{q_{n-1} q_{n-2}} + \cdots + \frac{p_1 q_0 - q_1 p_0}{q_1 q_0} + \frac{p_0}{q_0}$$
$$\overset{(A)}{=} \frac{(-1)^{n-1}}{q_n q_{n-1}} + \frac{(-1)^{n-2}}{q_{n-1} q_{n-2}} + \cdots + \frac{(-1)^0}{q_1 q_0} + a_0$$

where the sign theorem is applied in (A) and the last term $a_0 = p_0/q_0$ is the recursion theorem.

■

The next theorem is key for many estimations in the domain of continued fractions.

**Theorem 10** (*Monotony Theorem*)
Let $x_n \overset{\text{def}}{=} \frac{p_n}{q_n} = [a_0; a_1, \dots, a_n]$ denote the n-th convergent. Then:

$$x_{2n} < x_{2n+2}$$

and

$$x_{2n+1} > x_{2n+3}$$

I.e. even convergents are strictly monotonically increasing, and odd convergents are strictly monotonically decreasing.

Proof: We compute the following difference, where (A) uses again $-\frac{p_i}{q_i} + \frac{p_i}{q_i} = 0$:

$$
\begin{aligned}
x_n - x_{n-2} &= \frac{p_n}{q_n} - \frac{p_{n-2}}{q_{n-2}} \overset{(A)}{=} \frac{p_n}{q_n} - \frac{p_{n-1}}{q_{n-1}} + \frac{p_{n-1}}{q_{n-1}} - \frac{p_{n-2}}{q_{n-2}} \\
&= \frac{p_n q_{n-1} - q_n p_{n-1}}{q_n q_{n-1}} + \frac{p_{n-1} q_{n-2} - q_{n-1} p_{n-2}}{q_{n-1} q_{n-1}} \\
&\overset{(B)}{=} \frac{(-1)^{n-1}}{q_n q_{n-1}} + \frac{(-1)^{n-2}}{q_{n-1} q_{n-2}} = \frac{(-1)^{n-1} q_{n-2} + (-1)^{n-2} q_n}{q_n q_{n-1} q_{n-2}} \\
&= \frac{(-1)^{n-2} q_n - (-1)^{n-2} q_{n-2}}{q_n q_{n-1} q_{n-2}} = \frac{(-1)^{n-2}(q_n - q_{n-2})}{q_n q_{n-1} q_{n-2}} \\
&= \frac{(-1)^n (q_n - q_{n-2})}{q_n q_{n-1} q_{n-2}} \overset{(C)}{=} \frac{(-1)^n a_n q_{n-1}}{q_n q_{n-1} q_{n-2}} = \frac{(-1)^n a_n}{q_n q_{n-2}}
\end{aligned}
$$

(B) is because of the sign theorem, and (C) follows from $q_n = a_n q_{n-1} + q_{n-2}$, i.e. the recursion theorem.

Now, because of $a_n, q_n, q_{n-2} > 0$ it is $\frac{a_n}{q_n q_{n-2}} > 0$. Thus, $\frac{(-1)^n a_n}{q_n q_{n-2}} > 0$ for n even and $\frac{(-1)^n a_n}{q_n q_{n-2}} < 0$ for n odd. And this implies $x_n = \frac{(-1)^n a_n}{q_n q_{n-2}} + x_{n-2} > x_{n-2}$ for n even as well as $x_n = \frac{(-1)^n a_n}{q_n q_{n-2}} + x_{n-2} < x_{n-2}$ for n odd.

■

While even convergents are increasing and odd convergence are decreasing, all even convergents are smaller than all odd convergents. This is the content of the next very important theorem.

**Theorem 11** (*Convergents Comparison Theorem*)
For $0 \le 2n,\ 2m + 1 \le N$ it is $x_{2n} < x_{2m+1}$

Proof: As before, using the sign theorem in (A), we get

$$
x_n - x_{n-1} = \frac{p_n}{q_n} - \frac{p_{n-1}}{q_{n-1}} = \frac{p_n q_{n-1} - q_n p_{n-1}}{q_n q_{n-1}} \overset{(A)}{=} \frac{(-1)^{n-1}}{q_n q_{n-1}} = \frac{(-1)^{n-1}}{\beta_n}
$$

with $\beta_n := q_n q_{n-1}$. Because $q_n, q_{n-1} > 0$, it is $\beta_n > 0$, i.e. the sign of $\frac{(-1)^{n-1}}{\beta_n}$ is in fact $(-1)^{n-1}$.

Thus, $x_{2n+1} - x_{2n} = \frac{(-1)^{2n}}{\beta_{2n+1}} > 0$, and we get $x_{2n+1} = \frac{(-1)^{2n}}{\beta_{2n+1}} + x_{2n} > x_{2n}$ . This shows that an even convergent $x_{2n}$ is strictly smaller than its immediate succeeding odd convergent $x_{2n+1}$.

But what about an arbitrary odd convergent $x_{2m+1}$? For n < m the monotony theorem (theorem 11) yields $x_{2n} < x_{2m}$ and we showed before that $x_{2m} < x_{2m+1}$, thus, $x_{2n} < x_{2m+1}$.

For n > m, the monotony theorem yields $x_{2m+1} > x_{2n+1}$ and with $x_{2n+1} > x_{2n}$ we see $x_{2n} < x_{2m+1}$.

■

The following often used corollary computes the difference of two immediately succeeding convergents by mean of the denominators of the convergents, while the difference of the n-th convergent and the (n−2)-nd convergent adds the n-the element of the n-th convergent as a factor.

**Corollary 12**

$$\frac{p_n}{q_n} - \frac{p_{n-1}}{q_{n-1}} = \frac{(-1)^{n-1}}{q_n q_{n-1}} \tag{10}$$

and

$$\frac{p_n}{q_n} - \frac{p_{n-2}}{q_{n-2}} = \frac{(-1)^n a_n}{q_n q_{n-2}} \tag{11}$$

Proof: Equation (10) is the first equation from the proof of theorem 11. The second equation follows because of

$$\frac{p_n}{q_n} - \frac{p_{n-2}}{q_{n-2}} = \frac{p_n q_{n-2} - p_{n-2} q_n}{q_n q_{n-2}} \overset{(A)}{=} \frac{(-1)^n a_n}{q_n q_{n-2}},$$

where (A) is because of the second sign theorem 7.

■

We already saw that the even convergents are strictly monotonically increasing, that the odd convergents are strictly monotonically decreasing, and that each even convergent is less than all odd convergents. According to the next theorem the value of a continued fraction lies between the even convergents and the odd convergents, i.e. this value is larger than all even convergents and smaller than all odd convergents. The situation is depicted in figure 3.

Note, that the notion of the value of a continued fraction is defined by now for finite continued fractions. In section 4, this notion will also be defined for regular infinite continued fractions.

**Theorem 13** (*Nesting Theorem*)
Let x be the value of the continued fraction $[a_0; a_1, \ldots, a_N]$ and let $x_k$ be its convergents. Then:

$$\forall\, m, n < N : x_{2m} < x < x_{2n+1} \tag{12}$$

Proof: The value of x is the convergent with the highest index N, i.e. $x = x_N = [a_0; a_1, \ldots, a_N]$.

Let N = 2k be even. Since even convergents are strictly monotonically increasing, we know that $\forall 2m < N : x_{2m} < x_{2k} = x_N = x$, and according to the convergent comparison theorem 11 we know $\forall 2n+1 : x = x_N = x_{2k} < x_{2n+1}$.

Let N = 2k+1 be odd. Since odd convergents are strictly monotonically decreasing, we know that $\forall 2n+1 < N : x_{2n+1} > x_{2k+1} = x_N = x$, and according to the convergent comparison theorem 11 we know $\forall 2m : x = x_N = x_{2k+1} > x_{2m}$.

■

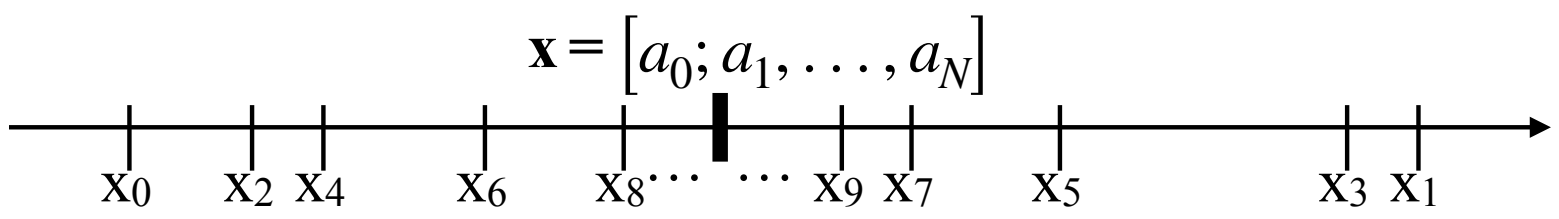

**Fig. 3**. Nesting of the value of a continued fraction by its convergents.

Because the value of a continued fraction is nested within its even convergents and odd convergents, the distance of this value from any of its convergents can be estimated by the distance of two consecutive convergents:

**Theorem 14** (*Distance Theorem*)
Let $x = [a_0; a_1, \ldots, a_N]$ and let $x_k$ be its convergents. Then:

$$\forall n : |x - x_n| < |x_{n-1} - x_n| \tag{13}$$

and

$$\forall n : |x - x_n| < |x_{n+1} - x_n| \tag{14}$$

Proof: Let n be even. Then $x_n < x < x_{n-1}$, i.e. $x - x_n < x_{n-1} - x_n$. Also, it is $x - x_n > 0$ and $x_{n-1} - x_n > 0$. Thus, $|x - x_n| < |x_{n-1} - x_n|$ for n even.

Now, let n be odd. It is $x_{n-1} < x < x_n$, which implies $x - x_n > x_{n-1} - x_n \Leftrightarrow -(x_n - x) > -(x_n - x_{n-1}) \Leftrightarrow x_n - x < x_n - x_{n-1}$. Because of $x_n - x > 0$ and $x_n - x_{n-1} > 0$, it is $|x_n - x| < |x_n - x_{n-1}| \Leftrightarrow |x - x_n| < |x_{n-1} - x_n|$ for n odd.

Together, this proves equation (13). Equation (14) is proven similarly.

■

Figure 3 shows the corresponding geometric situation for an even n.

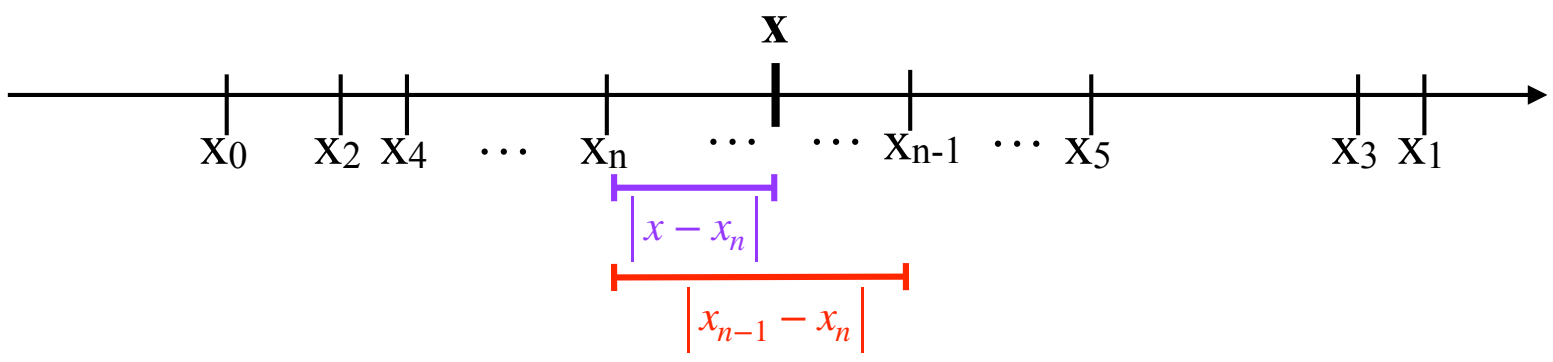

**Fig. 3**. The distance between two succeeding convergents is greater than the distance of a convergent and the value of its continued fraction.

Similarly, the difference between any two arbitrary convergents can be estimated by the difference of the convergents with the smaller index and its immediate predecessor:

**Theorem 15** (*Difference Theorem*)
Let $x = [a_0; a_1, \ldots, a_N]$ and let $x_k$ be its convergents. Then:

$$\forall\, m > n : \left|x_m - x_n\right| < \left|x_{n-1} - x_n\right| \tag{15}$$

Proof: Let n be even, e.g. n=2k.

Let m=2t be even. By theorem 11, even convergents are smaller than all odd convergents, i.e. $x_{2t} < x_{2k-1}$ for any $t \in \mathbb{N}$. Thus, $x_m - x_n = x_{2t} - x_{2k} < < x_{2k-1} - x_{2k} = x_{n-1} - x_n$.

Let m=2t−1 be odd. By the monotony theorem 10, odd convergents are strictly monotonically decreasing, i.e. $x_{2t-1} < x_{2k-1}$ for each t > k. Thus, $x_m - x_n = = x_{2t-1} - x_{2k} < x_{2k-1} - x_{2k} = x_{n-1} - x_n$.

For n odd, the proof is analogously.

■

The geometry of the last theorem is depicted in figure 4.

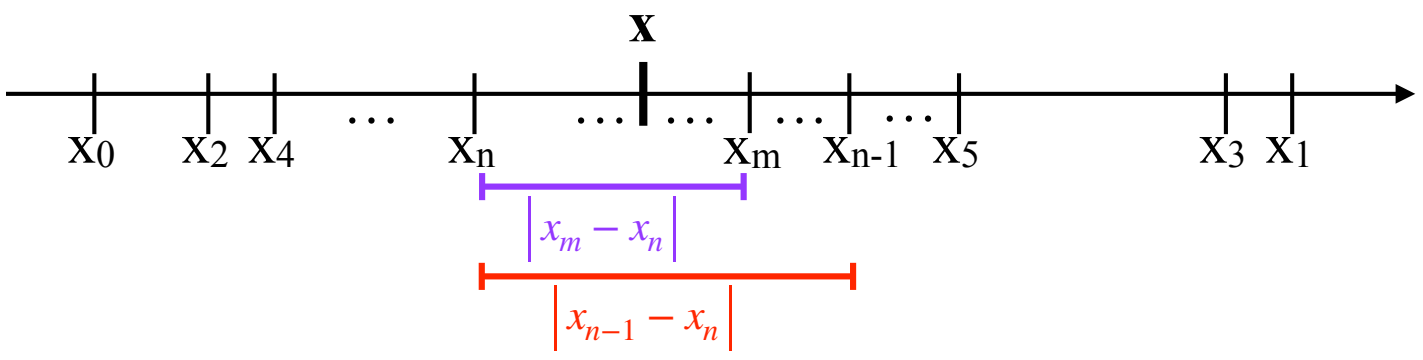

**Fig. 4**. The distance between any two convergents is smaller than the distance between the convergent with the smaller index and its immediate predecessor.

In several calculations the size of the denominator of a convergent must be estimated:

**Lemma 16** (*Size of Denominators*)
For the denominator $q_n$ of a convergent $\frac{p_n}{q_n} = [a_0; a_1, \ldots, a_n]$ the following holds:

$$\forall\, n : q_n \geq n \tag{16}$$

and

$$\forall\, n > 3 : q_n > n \tag{17}$$

Proof: By definition, $q_0 = 1 > 0$, and $q_1 = a_1 \geq 1$ because $a_i \in \mathbb{N}$, and finally,

$$q_2 \overset{(A)}{=} a_2 q_1 + q_0 \overset{(B)}{=} a_2 q_1 + 1 \overset{(C)}{\geq} q_1 + 1 \overset{(D)}{\geq} 2.$$

(A) holds because of the recursion theorem 4, (B) is by definition of $q_0$, (C) is because $a_2 \in \mathbb{N}$, and (D) has been seen just before (i.e. $q_1 \geq 1$). This proves the lemma for $n \leq 2$.

The proof for n ≥ 3 is by induction. It is

$$q_n \overset{(A)}{=} a_n q_{n-1} + q_{n-2} \overset{(B)}{\geq} q_{n-1} + q_{n-2} \overset{(C)}{\geq} q_{n-1} + (n-2) \overset{(D)}{\geq} q_{n-1} + 1 \overset{(E)}{\geq} n$$

where (A) is the recursion theorem, (B) is because of $a_n \in \mathbb{N}$, (C) is by induction hypothesis applied to $q_{n-2}$, (D) is because n ≥ 3, and (E) is by induction hypothesis applied to $q_{n-1}$. This proves equation (16).

Equation (17) is proven by induction again. Let n > 3. The argumentation is exactly as before, with the exception of (D):

$$q_n \overset{(A)}{=} a_n q_{n-1} + q_{n-2} \overset{(B)}{\geq} q_{n-1} + q_{n-2} \overset{(C)}{\geq} q_{n-1} + (n-2) \overset{(D)}{>} q_{n-1} + 1 \overset{(E)}{\geq} n$$

(D) holds because n > 3, i.e. $n - 2 > 1$.

■

In fact, denominators of a convergents grow much faster than the inequation $q_n > n$ may indicate:

**Lemma 17** (*Geometric Growth of Denominators*)
Let $q_n$ ($n \geq 2$) be the denominator of the convergent $\frac{p_n}{q_n} = [a_0; a_1, \ldots, a_n]$. Then:

$$q_n \geq 2^{\frac{n-1}{2}} \tag{18}$$

<u>Proof</u>: It is $q_k = a_k q_{k-1} + q_{k-2} > q_{k-1} + q_{k-2} \overset{(A)}{>} 2q_{k-2}$, with (A) because according to corollary 5, denominators are strictly monotonically increasing, i.e. $q_{k-1} > q_{k-2}$.

By induction, it is $q_{2k} \geq 2^k q_0$ and then $2^k q_0 \overset{(A)}{=} 2^k \overset{(B)}{\geq} 2^{\frac{(2k)-1}{2}}$ with (A) because $q_0 = 1$ and (B) follows from

$$2^k = 2^{\frac{2k}{2}} \geq \frac{1}{\sqrt{2}} 2^{\frac{2k}{2}} = 2^{\frac{2k}{2} - \frac{1}{2}} = 2^{\frac{2k-1}{2}}.$$

Similarly, by induction it is $q_{2k+1} \geq 2^k q_1$ and then $2^k q_1 \overset{(A)}{\geq} 2^k = 2^{\frac{(2k+1)-1}{2}}$ with (A) because of $q_1 \in \mathbb{N}$.

With $n = 2k$ and $n = 2k + 1$, respectively, follows equation (18).

■

## 4. Convergence of Infinite Regular Continuous Fractions

In section 2, we presented an algorithm to compute the continued fraction representation of a rational number. Next, we show how to compute such a representation for a non-rational number.

**Algorithm 18** (*Continued Fraction Representation of Non-Rational Numbers*)
Let $\alpha \in \mathbb{R}\backslash\mathbb{Q}$. Define:

- $\alpha_0 := \alpha$ and $b_0 := \lfloor \alpha_0 \rfloor$
- $\alpha_i := \dfrac{1}{\alpha_{i-1} - b_{i-1}}$ and $b_i := \lfloor \alpha_i \rfloor$ for i ≥ 1

Then, $\left[b_0; b_1, b_2, \ldots\right]$ is the continued fraction representation of α. Each $\alpha_i$ is called *i-th complete quotient* of α.

□

The above algorithm does not terminate, i.e. the continued fraction representation of a non-rational number is infinite. This is the content of the following

**Note 19**
In algorithm 18, it is $\alpha_i \notin \mathbb{Z}$

Proof (by induction):
n=0: Then by definition $\alpha_0 = \alpha \notin \mathbb{Z}$.
Induction hypothesis: $\alpha_n \notin \mathbb{Z}$
n → n+1: Assume $\alpha_n - b_n \in \mathbb{Z} \Rightarrow (\alpha_n - b_n) = k \in \mathbb{Z} \Rightarrow \alpha_n = k + b_n \in \mathbb{Z}$, which is a contradiction to the hypothesis! Thus, $\alpha_n - b_n \notin \mathbb{Z} \Rightarrow \alpha_{n+1} := \dfrac{1}{\alpha_n - b_n} \notin \mathbb{Z}$.

■

Figure 5 gives the computation of the continued fraction representation of $\sqrt{2}$:

$\sqrt{2} = 1{,}41421$

- $\alpha_0 = \sqrt{2}$ und $b_0 = \lfloor \sqrt{2} \rfloor = 1$
- $\alpha_1 = \dfrac{1}{\alpha_0 - b_0} = \dfrac{1}{0{,}41421} = 2{,}41421$ und $b_1 = 2$
- $\alpha_2 = \dfrac{1}{\alpha_1 - b_1} = \dfrac{1}{0{,}41421} = 2{,}41421$ und $b_2 = 2$

.
.
.

$\Rightarrow \sqrt{2} = [1; 2{,}2{,}2{,}...]$

**Fig. 5**. Computing the continued fraction of $\sqrt{2}$.

## 5. Bounds Expressed by Denominators of Convergents

In the following we give upper bounds and lower bounds of the approximations of a number by the convergents of its continued fraction representation by means of the denominators of the convergents.

First, we start with estimations of upper bounds:

**Lemma 20** (*Upper Bounds*)
Let $p_n/q_n$ be a convergent of the continued fraction representation of x. Then:

$$\left| x - \frac{p_n}{q_n} \right| < \frac{1}{q_n q_{n+1}} < \frac{1}{q_n^2} \le \frac{1}{n^2} \tag{19}$$

Proof: With $x_n = p_n/q_n$ it is $\left|x - x_n\right| < \left|x_{n+1} - x_n\right|$ (see theorem 14, equation (14)). According to corollary 12 (equation 10), it is

$$x_{n+1} - x_n = \frac{p_{n+1}}{q_{n+1}} - \frac{p_n}{q_n} = \frac{(-1)^n}{q_n q_{n+1}}.$$

Thus,

$$\left|x - x_n\right| < \left|x_{n+1} - x_n\right| = \left|\frac{(-1)^n}{q_n q_{n+1}}\right| = \frac{1}{q_n q_{n+1}} \overset{(A)}{<} \frac{1}{q_n^2} \overset{(B)}{\le} \frac{1}{n^2}$$

where (A) holds because of $q_{n+1} > q_n$ (corollary 5), and (B) is true because of $q_n \ge n$ (lemma 16).

■

An immediate consequence of this theorem is the convergence of the sequence of the convergents of a continued fraction to the value of the continued fraction. And this, by the way, is the origin of the name "convergents".

**Corollary 21**
The series $(p_n/q_n)$ of the convergents of the continued fraction representation of $x \in \mathbb{R}\backslash\mathbb{Q}$ converges to x:

$$\lim \frac{p_n}{q_n} = x$$

Proof: The claim follows immediately from $\left| x - \frac{p_n}{q_n} \right| < \frac{1}{n^2}$.

■

Often, two fractions are compared by means of their mediant ("mediant" means "somewhere in between").

**Definition 22**: For $a/b, c/d \in \mathbb{Q}$ and b, d > 0, the term $\frac{a+c}{b+d}$ is called the *mediant* of the two fractions. □

The following simple inequation is often used.

**Note 22** (*Mediant Property*)
Let $a/b, c/d \in \mathbb{Q}$ and b, d > 0 and $\frac{a}{b} < \frac{c}{d}$.
Then:

$$\frac{a}{b} < \frac{a+c}{b+d} < \frac{c}{d} \tag{20}$$

Proof: It is $\frac{a}{b} < \frac{c}{d} \Rightarrow ad < bc \Rightarrow bc - ad > 0$ and $b, d > 0 \Rightarrow b(b+d) > 0$.
This implies $\frac{a+c}{b+d} - \frac{a}{b} = \frac{b(a+c) - a(b+d)}{b(b+d)} = \frac{bc-ad}{b(b+d)} > 0$ and thus $\frac{a}{b} < \frac{a+c}{b+d}$. The inequation $\frac{a+c}{b+d} < \frac{c}{d}$ follows similarly.

■

Mediants of convergents that are weighted in a certain way are another important concept for computing bounds:

**Definition 23**: The term $x_{n,t} = \frac{tp_{n+1} + p_n}{tq_{n+1} + q_n}$ with $1 \le t \le a_{n+2}$ is called the (n,t)-th *semiconvergent*. □

Semiconvergents of an even n are strictly monotonically increasing, and semiconvergents of an odd n are strictly monotonically decreasing. This is the content of the following lemma.

**Lemma 23** (*Monotony of Semiconvergents*)
Let n be even; then $x_{n,t} < x_{n,t+1}$.
Let n be odd, the $x_{n,t} > x_{n,t+1}$

Proof: A simple calculation and the use of the sign theorem 6 results in

$$x_{n,t+1} - x_{n,t} = \frac{(t+1)p_{n+1} + p_n}{(t+1)q_{n+1} + q_n} - \frac{tp_{n+1} + p_n}{tq_{n+1} + q_n} = \frac{(-1)^n}{\left((t+1)q_{n+1} + q_n\right)\left(tq_{n+1} + q_n\right)}$$

The denominator of the last fraction of always positive. Thus, the last term is positive iff n is even (i.e. $x_{n,t+1} - x_{n,t} > 0$) , and it is negative is n is odd (i.e. $x_{n,t+1} - x_{n,t} < 0$).

■

In order to simplify proofs in what follows, the following conventions are used:

$$p_{-1} \stackrel{\text{def}}{=} 1 \text{ and } q_{-1} \stackrel{\text{def}}{=} 0 \tag{21}$$

With this, $x_{-1,1} = \frac{p_0 + p_{-1}}{q_0 + q_{-1}} = \frac{a_0 + 1}{1 + 0} = a_0 + 1$ becomes a semiconvergent. Now, $x_1 = \frac{p_1}{q_1} \overset{(A)}{=} \frac{a_1 a_0 + 1}{a_1} = a_0 + \frac{1}{a_1} \overset{(B)}{\leq} a_0 + 1 = x_{-1,1}$ where (A) is the recursion theorem and (B) follows because $a_1 \geq 1$, thus $x_1 \leq x_{-1,1}$.

Furthermore, it is $x_{-1,t} = \frac{t p_0 + p_{-1}}{t q_0 + q_{-1}} = \frac{t a_0 + 1}{t \cdot 1 + 0} = \frac{t a_0 + 1}{t} = a_0 + \frac{1}{t}$ for $1 \leq t \leq a_1$.

Putting things together, it is

$$x_{-1,1} = a_0 + 1 > a_0 + \frac{1}{2} > \cdots > a_0 + \frac{1}{a_1} = x_1 \tag{22}$$

Based on this we can refine figure 3, which depicts the nesting and ordering of convergents by including semiconvergents: Between two succeeding convergents (e.g. $x_n$ and $x_{n+2}$ in figure 6) the corresponding semiconvergents ordered according to lemma 23 are nested (in increasing order as shown for an even n in figure 6). Furthermore, beyond $x_1 = a_0 + \frac{1}{a_1}$, the semiconvergents $x_{-1,t}$ are added.

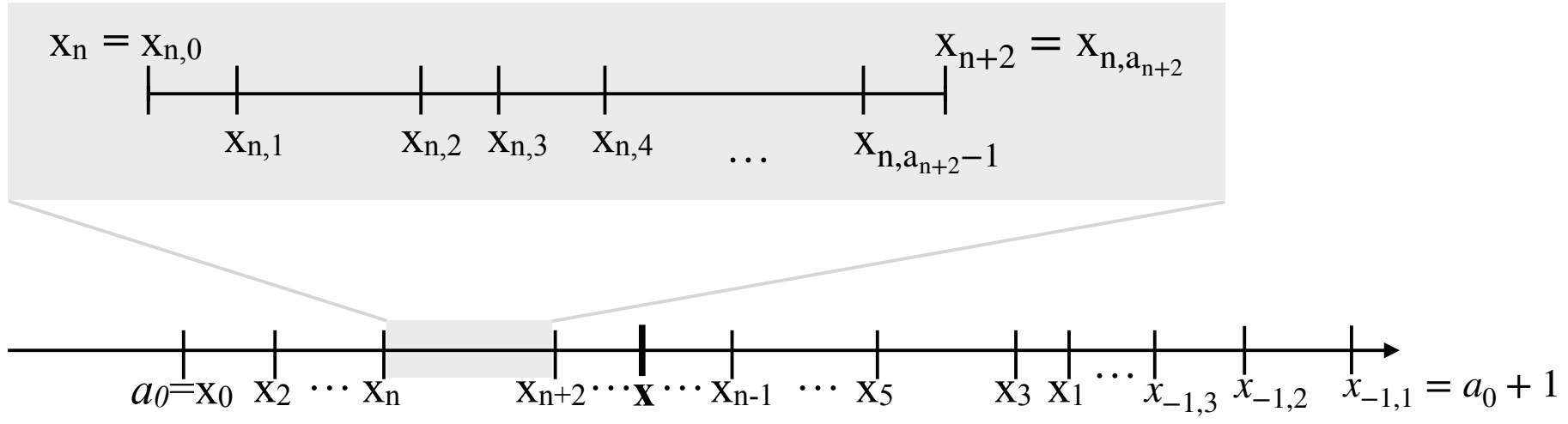

**Fig. 6**. Nesting of convergents and semiconvergents (n even).

Now, we are prepared to prove a lower bound of the approximation of a number by the convergents of its continued fraction representation by means of the denominators of the convergents.

**Lemma 24** (*Lower Bound*)
Let $p_n/q_n$ be a convergent of the continued fraction representation of x. Then:

$$\left| x - \frac{p_n}{q_n} \right| > \frac{1}{(q_n + q_{n+1})\, q_n} \tag{23}$$

Proof: The proof is based on the following claims.

Claim 1: n even $\Rightarrow \frac{p_n}{q_n} < \frac{p_{n+1} + p_n}{q_{n+1} + q_n} < x < \frac{p_{n+1}}{q_{n+1}}$

Proof: $\frac{p_{n+1}+p_n}{q_{n+1}+q_n}$ is the mediant of $\frac{p_{n+1}}{q_{n+1}}$ and $\frac{p_n}{q_n}$. Thus, the mediant property (note 22) shows that $\frac{p_n}{q_n} < \frac{p_{n+1}+p_n}{q_{n+1}+q_n} < \frac{p_{n+1}}{q_{n+1}}$. Then:

$$\frac{p_n}{q_n} < \frac{p_{n+1}+p_n}{q_{n+1}+q_n} \overset{(A)}{<} \frac{2p_{n+1}+p_n}{2q_{n+1}+q_n} < \cdots < \frac{a_{n+2}p_{n+1}+p_n}{a_{n+2}q_{n+1}+q_n} \overset{(B)}{=} \frac{p_{n+2}}{q_{n+2}},$$

where (A) follows by the monotony of even semiconvergents (lemma 23), and (B) is the recursion theorem. Because of the nesting theorem 13 (note, that n+2 is even and n+1 is odd) it is $\frac{p_{n+2}}{q_{n+2}} < x < \frac{p_{n+1}}{q_{n+1}}$. This proves claim 1. $\square_{\text{(claim 1)}}$

Claim 2: n odd $\Rightarrow \frac{p_{n-1}}{q_{n-1}} < x < \frac{p_{n+1}+p_n}{q_{n+1}+q_n} < \frac{p_n}{q_n}$

Proof: As before, $\frac{p_{n+1}+p_n}{q_{n+1}+q_n}$ is the mediant of $\frac{p_{n+1}}{q_{n+1}}$ and $\frac{p_n}{q_n}$. Because n is odd, it is $\frac{p_{n+1}}{q_{n+1}} < \frac{p_n}{q_n}$ (nesting theorem 13). Thus, $\frac{p_{n+1}}{q_{n+1}} < \frac{p_{n+1}+p_n}{q_{n+1}+q_n} < \frac{p_n}{q_n}$ because of the mediant property (note 22). Then:

$$\frac{p_n}{q_n} > \frac{p_{n+1}+p_n}{q_{n+1}+q_n} \overset{(A)}{>} \frac{2p_{n+1}+p_n}{2q_{n+1}+q_n} > \cdots > \frac{a_{n+2}p_{n+1}+p_n}{a_{n+2}q_{n+1}+q_n} \overset{(B)}{=} \frac{p_{n+2}}{q_{n+2}},$$

where (A) follows by the monotony of odd semiconvergents (lemma 23), and (B) is the recursion theorem. Because of the nesting theorem 13 (note, that n−1 is even and n+2 is odd) it is $\frac{p_{n-1}}{q_{n-1}} < x < \frac{p_{n+2}}{q_{n+2}}$, and because n is odd, it is $\frac{p_{n+2}}{q_{n+2}} < \frac{p_n}{q_n}$. This proves claim 2. $\square_{\text{(claim 2)}}$

With claim 1, for even n it is $\frac{p_n}{q_n} < \frac{p_{n+1}+p_n}{q_{n+1}+q_n} < x \Rightarrow x - \frac{p_n}{q_n} > \frac{p_n+p_{n+1}}{q_n+q_{n+1}} - \frac{p_n}{q_n}$.

For n odd and claim 2 it is $x < \frac{p_{n+1}+p_n}{q_{n+1}+q_n} < \frac{p_n}{q_n} \Rightarrow \frac{p_n}{q_n} - x > \frac{p_n}{q_n} - \frac{p_n+p_{n+1}}{q_n+q_{n+1}} \Leftrightarrow$ $-\left(x - \frac{p_n}{q_n}\right) > -\left(\frac{p_n+p_{n+1}}{q_n+q_{n+1}} - \frac{p_n}{q_n}\right)$.

Thus, for any $k \in \mathbb{N}$: $\left| x - \frac{p_k}{q_k} \right| > \left| \frac{p_{k+1}+p_k}{q_{k+1}+q_k} - \frac{p_k}{q_k} \right|$. Next, we compute

$$\frac{p_k+p_{k+1}}{q_k+q_{k+1}} - \frac{p_k}{q_k} = \frac{(p_k+p_{k+1})\,q_k - (q_k+q_{k+1})\,p_k}{(q_k+q_{k+1})\,q_k} = \frac{p_{k+1}q_k - p_kq_{k+1}}{(q_k+q_{k+1})\,q_k}$$

$$\overset{(A)}{=} \frac{(-1)^k}{(q_k+q_{k+1})\,q_k}$$

where (A) is the sign theorem 6.

This implies $\left| x - \frac{p_k}{q_k} \right| > \left| \frac{(-1)^k}{(q_k + q_{k+1})\, q_k} \right| = \frac{1}{(q_k + q_{k+1})\, q_k}$.

■

Because of $q_{k+1} > q_k$ (corollary 5) it is

$$q_k + q_{k+1} < 2q_{k+1} \Leftrightarrow \frac{1}{2q_{k+1}} < \frac{1}{q_k + q_{k+1}} \Leftrightarrow \frac{1}{2q_k q_{k+1}} < \frac{1}{(q_k + q_{k+1})\, q_k}.$$

Using the last inequality in lemma 24 (Lower Bound) and using lemma 20 (Upper bounds) we get the concluding theorem of this section:

In summary, we have proved the following

**Theorem 25** (*Bounds of Approximations by Convergents*)
Let $p_k / q_k$ be a convergent of the continued fraction representation of x. Then:

$$\frac{1}{2q_k q_{k+1}} < \frac{1}{(q_k + q_{k+1})\, q_k} < \left| x - \frac{p_k}{q_k} \right| < \frac{1}{q_k q_{k+1}} < \frac{1}{q_k^2} \qquad (24)$$

■

## 6. Best Approximations

Our goal is to approximate a real number by a rational number as good as possible while keeping the denominator of the rational number "small". Keeping the denominator small is important because in practice every real number can only be given up to a certain degree of precision, and this is achieved by means of a huge denominator and corresponding numerator. I.e. approximating a real number by a rational number with huge denominator is canonical, but finding a small denominator is a problem.

This is captured by the following

**Definition 26**: A fraction $p/q \in \mathbb{Q}$ is called *best approximation* (*of the first kind*) of $\alpha \in \mathbb{R} :\Leftrightarrow \forall\, c/d \in \mathbb{Q} : d \le q \Rightarrow \left| \alpha - \frac{c}{d} \right| > \left| \alpha - \frac{p}{q} \right|$ (assuming c/d ≠ p/q). □

Often, the addition "of the first kind" is omitted. By definition, a best approximation of a real number can only be improved if the denominator of the given approximation is increased.

If p/q is a best approximation of α then $\left| \alpha - \frac{p}{q} \right| = \frac{1}{q} \left| q\alpha - p \right|$ is small and, thus, $\left| q\alpha - p \right|$ is small. Measuring the goodness of an approximation this way results in the next

**Definition 27**: A fraction $p/q \in \mathbb{Q}$ is called *best approximation of the second kind* of $\alpha \in \mathbb{R} :\Leftrightarrow \forall\, c/d \in \mathbb{Q} : d \leq q \Rightarrow |d\alpha - c| > |q\alpha - p|$ (assuming c/d ≠ p/q). □

The question is whether every best approximation is also a best approximation of the second kind. Now, 1/3 is best approximation of 1/5 because the only possible fractions for c/d with d ≤ 3 = q, are 0, 1/2, 2/3 and 1, and these numbers satisfy $\left|\frac{1}{5} - \frac{c}{d}\right| > \left|\frac{1}{5} - \frac{1}{3}\right|$.

Next we observe that $\left|1 \cdot \frac{1}{5} - 0\right| < \left|3 \cdot \frac{1}{5} - 1\right|$ with $1 < 3$. Thus, with $d = 1$ and $q = 3$ (i.e. $d < q$) and $\alpha = 1/5$, we found a fraction $c/d = 0/1$ with $|d\alpha - c| < |q\alpha - p|$! As a consequence, although 1/3 is a best approximations of the first kind of 1/5 it is not a best approximation of the second kind.

Thus, not all best approximations of the first kind are best approximations of the second kind. But the reverse holds true:

**Lemma 28** (*Every 2nd Kind Best Approximation is a 1st Kind Best Approximation*)
If $p/q \in \mathbb{Q}$ is a best approximation of the second kind of $\alpha \in \mathbb{R}$, then $p/q$ is also a best approximation of the first kind of α.

Proof (by contradiction): Assume p/q is not best approximation of the first kind. Then, $\left|\alpha - \frac{c}{d}\right| \leq \left|\alpha - \frac{p}{q}\right|$ for a fraction c/d with d < q. Multiplying both inequations results in $d\left|\alpha - \frac{c}{d}\right| \leq q\left|\alpha - \frac{p}{q}\right| \Leftrightarrow |d\alpha - c| \leq |q\alpha - p|$, which is a contradiction because p/q is a best approximation of the second kind.

■

The next simple estimation about the distance of two fractions by means of the product of their denominators is often used.

**Note 29** (*Distance of Fractions*)
Let $\frac{a}{b}, \frac{p}{q} \in \mathbb{Q}$ with $\frac{a}{b} \neq \frac{p}{q}$. Then:

$$\left|\frac{p}{q} - \frac{a}{b}\right| \geq \frac{1}{qb} \tag{25}$$

Proof: With $a, p \in \mathbb{Z}$ and $b, q \in \mathbb{N}$ it is $pb - aq \in \mathbb{Z}$. Also, $pb - aq \neq 0$ because otherwise $pb = aq \Leftrightarrow \frac{p}{q} = \frac{a}{b}$ which contradicts the premise. Thus, $|pb - aq| \in \mathbb{N}$, i.e. $|pb - aq| \geq 1$. This implies

$$\left|\frac{p}{q}-\frac{a}{b}\right| = \left|\frac{pb-aq}{qb}\right| = \frac{|pb-aq|}{|qb|} \geq \frac{1}{qb},$$

where $|qb| = qb$ because $b, q \in \mathbb{N}$.

■

Next we prove that every best approximation of the second kind is a convergent.

**Theorem 30** (*2nd Kind Best Approximations are Convergents*)
Let $a/b$ be a best approximation of the second kind of $x \in \mathbb{R}$, and let $x = [a_0; a_1, \cdots]$ be the continued fraction representation of x.

Then $a/b$ is a convergent of x.

Proof: Being a best approximation of the second kind of x, $a/b$ satisfies by definition $|dx - c| > |bx - a|$ for d ≤ b.

Claim 1: $\frac{a}{b} \geq a_0 = x_0$

Proof (by contradiction): Assume $\frac{a}{b} < a_0 \Rightarrow -a_0 < -\frac{a}{b} \Rightarrow x - a_0 < x - \frac{a}{b}$, thus $|x - a_0| < \left|x - \frac{a}{b}\right| \overset{(A)}{\leq} b\left|x - \frac{a}{b}\right| = |bx - a|$, where (A) holds because $b \in \mathbb{N}$, i.e. 1 ≤ b. This implies $|1 \cdot x - a_0| \leq |bx - a|$, which contradicts $|dx - c| > |bx - a|$ for d ≤ b (with $d = 1 \leq b$ and $c = a_0$). This means that $\frac{a}{b} \geq a_0 = \frac{a_0}{1} \overset{(B)}{=} \frac{q_0}{q_0} = x_0$, (B) is because of the recursion theorem. $\square_{\text{(claim 1)}}$

Thus, the geometric situation is as depicted in figure 7, i.e. $a/b$ is in the grey shaded area being greater than or equal to the convergent $x_0$. This will be refined in what follows.

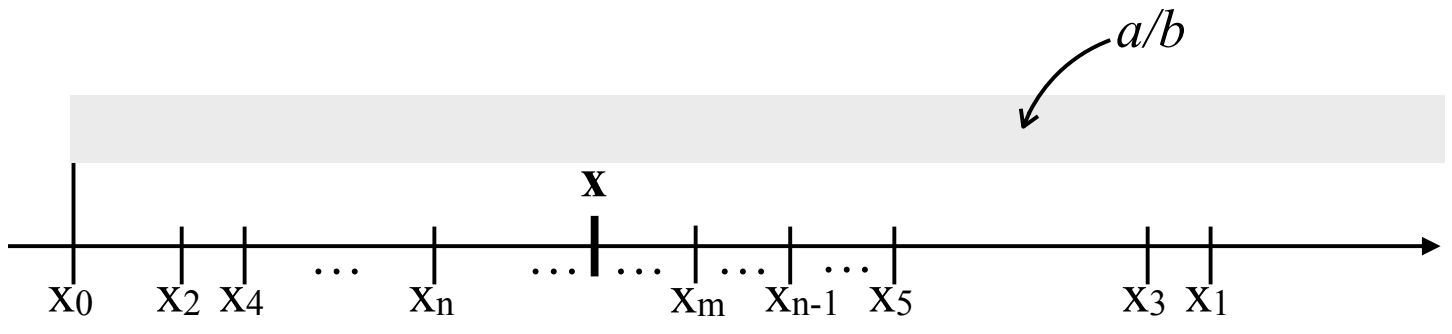

**Fig. 7**. Any best approximation of the second kind is in the grey shaded area, i.e. greater than or equal the convergent $x_0$.

Next, we proceed with a proof by contradiction assuming that $a/b$ is no convergent of x.

Assumption: $\frac{a}{b} \neq \frac{q_k}{q_k} = x_k$ for $k \in \mathbb{N}$

According to claim 1, $\frac{a}{b} \geq a_0 = x_0$. Thus, one of the following must hold:

$$\frac{a}{b} \in \left]\frac{p_{k-1}}{q_{k-1}}, \frac{p_{k+1}}{q_{k+1}}\right[ \text{ for k} \geq 1 \qquad (*)$$

or

$$\frac{a}{b} > \frac{p_1}{q_1} = x_1 \qquad (**)$$

This situation is shown in figure 8.

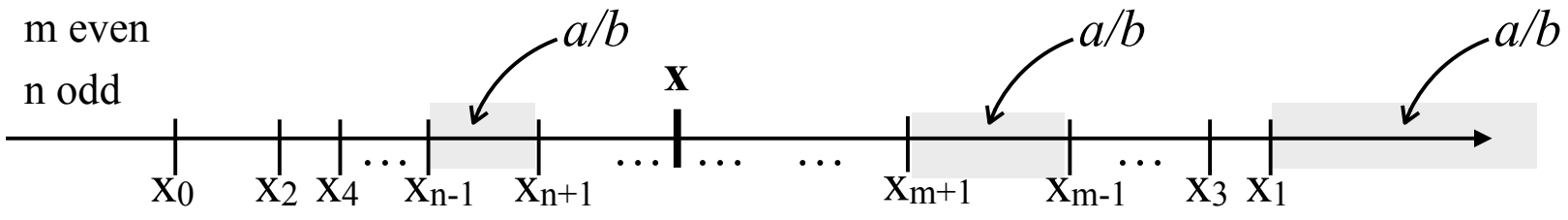

**Fig. 8**. If a best approximation of the second kind is not a convergent, it is within the indicated grey shaded areas.

<u>Case (*)</u>: If (*) is true, then

$$\left|\frac{a}{b} - \frac{p_{k-1}}{q_{k-1}}\right| < \left|x - \frac{p_{k-1}}{q_{k-1}}\right| \overset{(*)}{<} \left|\frac{p_k}{q_k} - \frac{p_{k-1}}{q_{k-1}}\right| \overset{(C)}{=} \frac{1}{q_k q_{k-1}},$$

where (*) is theorem 14, Eq. (14), and (C) is from corollary 12, Eq. (10). Furthermore, $\left|\frac{a}{b} - \frac{p_{k-1}}{q_{k-1}}\right| \overset{(D)}{\geq} \frac{1}{b q_{k-1}}$, with (D) because of note 29 (Distance of Fractions).

Together, $\frac{1}{b q_{k-1}} \leq \left|\frac{a}{b} - \frac{p_{k-1}}{q_{k-1}}\right| < \frac{1}{q_k q_{k-1}} \Rightarrow \frac{1}{b} < \frac{1}{q_k} \Rightarrow b > q_k$ (§).

Also, if (*) is true, then $\left|x - \frac{a}{b}\right| \geq \left|\frac{p_{k+1}}{q_{k+1}} - \frac{a}{b}\right| \overset{(E)}{\geq} \frac{1}{b q_{k+1}}$, where (E) is again using note 29. This implies $b\left|x - \frac{a}{b}\right| \geq \frac{1}{q_{k+1}} \Rightarrow \left|bx - a\right| \geq \frac{1}{q_{k+1}}$ (§§).

Lemma 20 (Upper Bounds) tells us that $\left|x - \frac{p_k}{q_k}\right| < \frac{1}{q_k q_{k+1}}$ which is equivalent to $q_k\left|x - \frac{p_k}{q_k}\right| < \frac{1}{q_{k+1}} \Leftrightarrow \left|q_k x - p_k\right| < \frac{1}{q_{k+1}} \Rightarrow \left|q_k x - p_k\right| < \left|bx - a\right|$ (see (§§) just before). Since, $q_k < b$ (see (§) above), this is a contradiction to $a/b$ being a best approximation of the second kind of x. Thus, case (*) does not occur.

<u>Case (**)</u>: This case is shown in figure 9. Then, $\left|x - \frac{a}{b}\right| > \left|\frac{p_1}{q_1} - \frac{a}{b}\right| \overset{(F)}{=} \frac{1}{b q_1}$, where (F) again uses note 29. This implies $\left|bx - a\right| > \frac{1}{q_1} \overset{(G)}{=} \frac{1}{a_1}$ (§§§) with (G) using the recursion theorem.

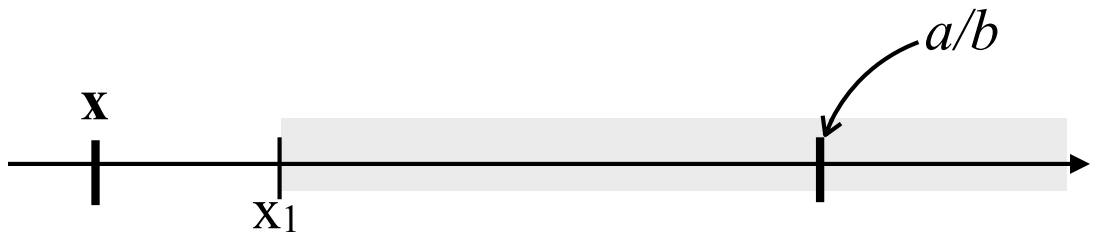

**Fig. 9**. Pictorial representation of case (**).

Now, $x - a_0 = \cfrac{1}{a_1 + \cfrac{1}{a_2 + \ddots}} \leq \frac{1}{a_1}$, where the last inequality holds because of $\frac{1}{a_2 + \ddots} > 0$, thus $\left|x - a_0\right| \leq \frac{1}{a_1} \overset{(H)}{<} \left|bx - a\right|$, (H) based on (§§§) before. This means that $\left|1 \cdot x - a_0\right| < \left|bx - a\right|$ with $1 \leq b$, i.e. $a/b$ is no best approximation of the second kind of x - which is a contradiction. Thus, case (**) does not occur either.

Consequently, the assumption is wrong and there is a $k \in \mathbb{N}$ with $\frac{a}{b} = \frac{q_k}{q_k} = x_k$, i.e. $a/b$ is a convergent.

■

So, every best approximation of the second kind is a convergent. The next theorem proves the reverse, i.e. that every convergent is a best approximation of the second kind.

**Theorem 31** (*Lagrange, 1798 — Convergents are 2nd Kind Best Approximations*) Let $p_n/q_n$ be a convergent of $x = \left[a_0; a_1, \cdots, a_N\right]$, $x \neq a_0 + \frac{1}{2}$ and $n \neq 0$. Then, for $d \leq q_n$ and $\frac{c}{d} \neq \frac{p_n}{q_n}$ it is $\left|dx - c\right| > \left|q_n x - p_n\right|$, i.e. the convergent is a best approximation of the second kind of x.

The case $x = a_0 + \frac{1}{2}$ and $n = 0$ is excluded because the convergent $\frac{p_0}{q_0} = \frac{a_0}{1}$ is not a best approximation of the second kind of $x = a_0 + \frac{1}{2}$: it is $\left|1 \cdot x - \left(a_0 + 1\right)\right| = \left|a_0 + \frac{1}{2} - a_0 - 1\right| = \frac{1}{2}$ and $\left|1 \cdot x - a_0\right| = \left|a_0 + \frac{1}{2} - a_0\right| = \frac{1}{2}$, which implies $\left|1 \cdot x - \left(a_0 + 1\right)\right| = \left|1 \cdot x - a_0\right|$. Setting $d := 1 \leq q_0$, $c := a_0 + 1$ results in $\left|d \cdot x - c\right| = \left|1 \cdot x - \left(a_0 + 1\right)\right| = \left|1 \cdot x - a_0\right| = \left|q_0 \cdot x - p_0\right|$. If $\frac{p_0}{q_0}$ would be a best approximation of the second kind of x, then $\left|1 \cdot x - \left(a_0 + 1\right)\right| > \left|1 \cdot x - a_0\right|$ would hold.

The proof of Lagrange's theorem is very technical. First, the expression $\left|y_0 x - z_0\right|$ is analyzed to find the smallest integral numbers $y_0$ and $z_0$ such that the expression is minimized under the constraint $y_0 \in \left\{q_0, \ldots, q_k\right\}$, i.e. $y_0$ is a denominator of a convergent. It is shown both, that $z_0/y_0$ is a best approximation of the second kind of x, and finally that $z_0 = p_k$ and $y_0 = q_k$.

Proof: Let $k \in \mathbb{Z}$ and $p_k/q_k$ be a convergent. First, we are looking for the smallest numbers $y_0, z_0 \in \mathbb{Z}$ with $y_0 \in \{q_0, \dots, q_k\}$ such that $\left|y_0 x - z_0\right|$ is minimal.

Step 1: Pick an arbitrary $z \in \mathbb{Z}$, and based on this we determine $y_0 \in \{q_0, \dots, q_k\}$.

It is $\min_y \left|yx - z\right| = 0 \Leftrightarrow y = \frac{z}{x}$, but in general $y \notin \mathbb{Z}$. Looking for a solution $y_0 \in \{q_0, \dots, q_k\} \subseteq \mathbb{Z}$ that minimizes $\left|y_0 x - z\right|$ results in the following potential positions of $z/x$ with respect to the denominators $q_0, \dots, q_k$ (see figure 10):

- Case 1: $z/x > q_k$. Then, $y_0 = q_k$ is the solution.
- Case 2: $z/x < q_0$. Then, $y_0 = q_0$ is the solution.

Let $q_i \leq z/x \leq q_{i+1}$ for $1 \leq i \leq k$.

- Case 3: For $\left|q_{i+1}x - z\right| < \left|q_i x - z\right|$ (i.e. $z/x$ is closer to $q_{i+1}$ than to $q_i$), $y_0 = q_{i+1}$ is the solution, and for $\left|q_{i+1}x - z\right| > \left|q_i x - z\right|$ (i.e. $z/x$ is closer to $q_i$ than to $q_{i+1}$), $y_0 = q_i$ is the solution.
- Case 4: For $\left|q_{i+1}x - z\right| = \left|q_i x - z\right|$ (i.e. $z/x$ is exactly in the middle between $q_i$ and $q_{i+1}$), $y_0 = q_i$ is the solution because $q_i < q_{i+1}$, and we are looking for the smallest $y_0$.

Especially, $y_0 \geq q_0 = 1$. $\square_{\text{(step 1)}}$

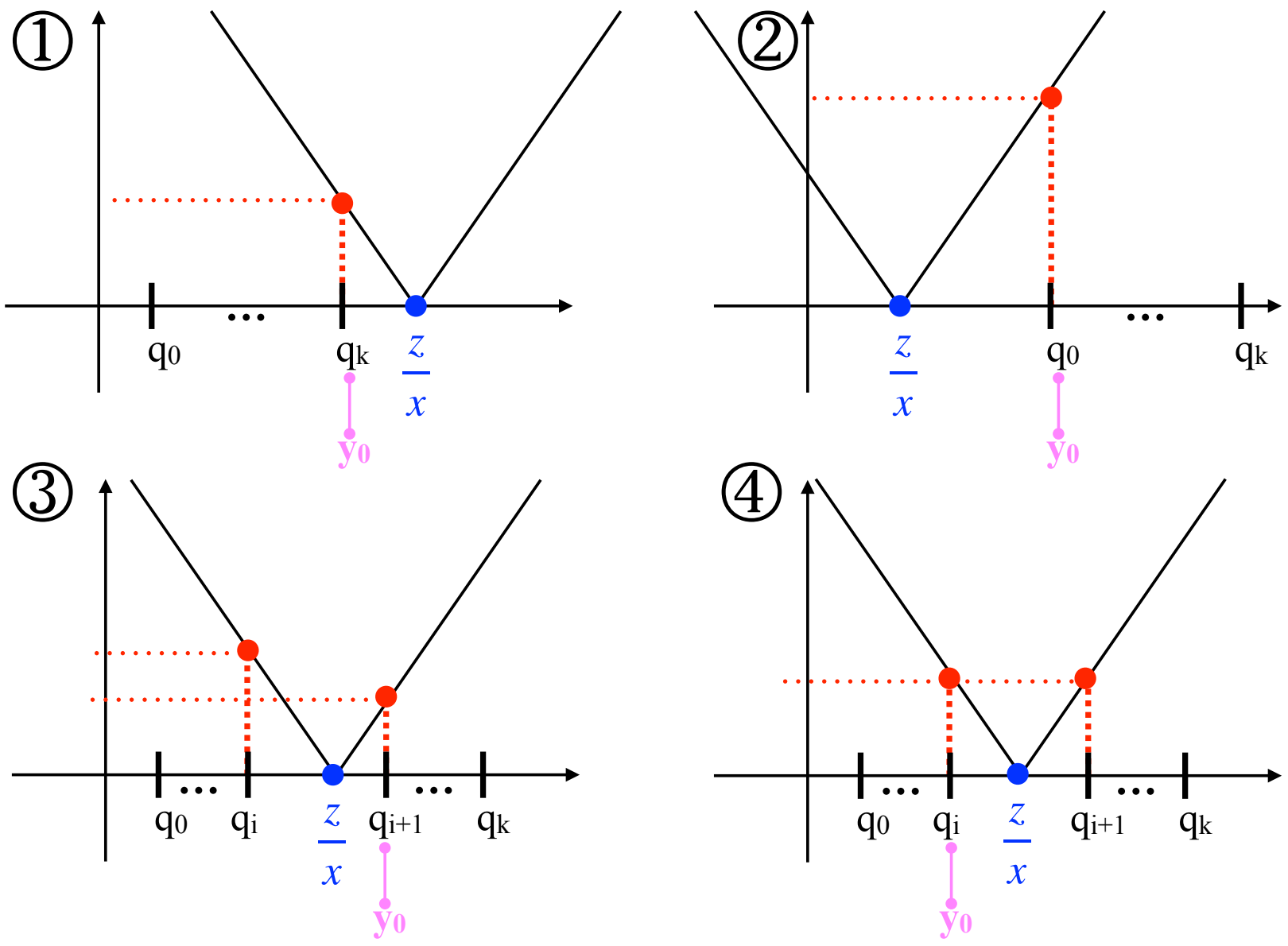

**Fig. 10**. The potential positions of $z/x$ with respect to the denominators $q_0, \dots, q_k$.

Step 2: Based on the $y_0$ found, we determine $z_0$ next. It is $\min_z \left| y_0 x - z \right| = 0$ $\Leftrightarrow z = y_0 x$, but in general $z \notin \mathbb{Z}$. In solving the minimization problem within $\mathbb{Z}$ (i.e. $z_0 := \operatorname{argmin}_{z \in \mathbb{Z}} \left| y_0 x - z \right|$), the following cases can be distinguished (see figure 11):

- Case 0: It may happen that $y_0 x \in \mathbb{Z}$. Then, choose $z_0 = y_0 x$ .
- Case 1: $y_0 x$ is between two integral numbers s and t, i.e. $s < y_0 x < t$ . For $\left| y_0 x - s \right| > \left| y_0 x - t \right|$ (i.e. $y_0 x$ is closer to t than to s), $z_0 = t$ is the solution; and for $\left| y_0 x - s \right| < \left| y_0 x - t \right|$ (i.e. $y_0 x$ is closer to s than to t), $z_0 = s$ is the solution.
- Case 2: For $\left| y_0 x - s \right| = \left| y_0 x - t \right|$ (i.e. $y_0 x$ is exactly in the middle between t and s), $z_0 = s$ is the solution because $s < t$, and we are looking for the smallest $z_0$. $\square$(step 2)

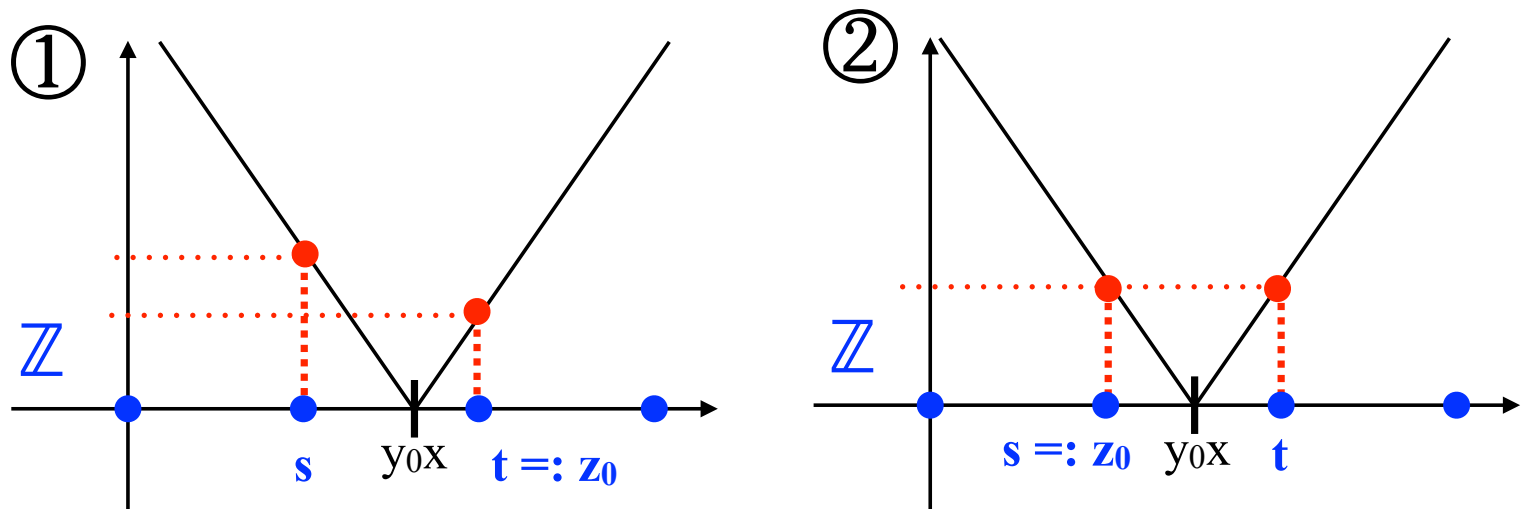

**Fig. 11**. The potential positions of $y_0 x$.

Claim 1: $z_0$ is uniquely determined.

Proof (by contradiction): Assume there exists a $\tilde{z} \in \mathbb{Z}$ with $\tilde{z}_0 \neq z_0$ and $\left| x - \frac{z_0}{y_0} \right| = \left| x - \frac{\tilde{z}_0}{y_0} \right|$. This can only happen iff one term is positive and the other is negative, i.e. for example if $x - \frac{z_0}{y_0} > 0$ and $x - \frac{\tilde{z}_0}{y_0} < 0$, and then $x - \frac{z_0}{y_0} = \frac{\tilde{z}_0}{y_0} - x$, i.e. $x = \frac{z_0 + \tilde{z}_0}{2y_0}$.

As an intermediate step we prove

Claim 2: $z_0 + \tilde{z}_0$ and $2y_0$ are co-prime, i.e. $\gcd\left(z_0 + \tilde{z}_0, 2y_0\right) = 1$

Proof (by contradiction): Let $\tilde{z}_0 + z_0 = Lp$ and $2y_0 = Lq$ with $L > 1$. Then: $x = \frac{z_0 + \tilde{z}_0}{2y_0} = \frac{Lp}{Lq} \Rightarrow x = \frac{p}{q}$ and thus

$$\left| qx - p \right| = \left| q\frac{p}{q} - p \right| = 0 \, (\S).$$

Assume $L > 2$. Then, with $2y_0 = Lq$ and $L/2 > 1$ it follows

$$y_0 = \frac{L}{2}q > q\ (\S\S).$$

Now, $y_0$ has been determined in step 1 to satisfy $y_0 = \operatorname{argmin}_y \left|yx - z\right|$ for a given z, especially for $z = p$, i.e. $y_0 = \operatorname{argmin}_y \left|yx - p\right|$. Because $0 = \min_y \left|yx - p\right|$ and $\left|qx - p\right| = 0$, it must be $q = y_0$: contradiction because $q < y_0$ according to (§§) before. Thus, $1 < L \leq 2$, i.e. $L = 2$.

With $L = 2$ and $2y_0 = Lq$ we get $y_0 = q$, which implies by definition of $z_0$: $\left|qx - p\right| = \left|y_0x - p\right| > \left|y_0x - z_0\right|$. But $\left|qx - p\right| = 0$ (see (§) above), thus $0 > \left|y_0x - z_0\right|$, which is a contraction. $\square_{\text{(claim 2)}}$

We continue the proof of claim 1: It is $\frac{z_0 + \tilde{z}_0}{2y_0} = x$ and also $x = \frac{p_N}{q_N}$, i.e. $\frac{z_0 + \tilde{z}_0}{2y_0} = \frac{p_N}{q_N}$. Because $\gcd\left(z_0 + \tilde{z}_0, 2y_0\right) = 1$ according to claim 2, it follows that $p_N = z_0 + \tilde{z}_0$ and $q_N = 2y_0$.

Now, let $N \geq 2$. Then it is $2y_0 = q_N \overset{(A)}{=} a_Nq_{N-1} + q_{N-2}$ ((A) uses the recursion theorem 4), and with note 2 is is $a_N \geq 2$. Thus, $2y_0 \geq 2q_{N-1} + q_{N-2}$ $\Rightarrow y_0 \geq q_{N-1} + \frac{q_{N-2}}{2} \Rightarrow q_{N-1} \leq y_0 - \frac{q_{N-2}}{2} \overset{(B)}{<} y_0$ ((B) is because of $q_{N-2} > 0$). Now:

$$\left|q_{N-1}x - p_{N-1}\right| = \left|q_{N-1}\frac{p_N}{q_N} - p_{N-1}\right| = \frac{1}{q_N}\left|q_{N-1}p_N - p_{N-1}q_N\right| \overset{(C)}{=} \frac{1}{q_N} = \frac{1}{2y_0} \overset{(D)}{\leq} \frac{1}{2}$$

where (C) holds because of the sign theorem and (D) because $y_0 \geq 1$ (see at the end of the proof of step 1).

Furthermore,

$$\left|y_0x - z_0\right| = \left|y_0\frac{z_0 + \tilde{z}_0}{2y_0} - z_0\right| = \left|\frac{z_0 + \tilde{z}_0}{2} - z_0\right| = \frac{1}{2}\left|z_0 + \tilde{z}_0 - 2z_0\right|$$

$$= \frac{1}{2}\left|\tilde{z}_0 - z_0\right| \overset{(E)}{\geq} \frac{1}{2} \qquad (\S\S\S)$$

where (E) is true because $\tilde{z}_0 \neq z_0$ and thus $\left|\tilde{z}_0 - z_0\right| \geq 1$ for integral numbers $\tilde{z}_0$ and $z_0$. Together, we got $\left|y_0x - z_0\right| \geq \frac{1}{2} \geq \left|q_{N-1}x - p_{N-1}\right|$, which is a contradiction to the choice of $y_0$ and $z_0$! This proves claim 1 for $N \geq 2$.

Now, let $N = 1$ and choose $a_1 = 2$ (based on note 2, the highest element of a continued fraction is always greater than or equal 2, thus $a_1 \geq 2$). Then

$$x = \left[a_0; a_1\right] = \frac{p_1}{q_1} \overset{(F)}{=} \frac{a_1a_0 + 1}{a_1} = \frac{2a_0 + 1}{2} = a_0 + \frac{1}{2}$$

((F) is the recursion theorem) which has been excluded from the theorem.

Thus, let $N = 1$ and $a_1 > 2$. Then

$$\left|1 \cdot x - a_0\right| \overset{(G)}{=} \left|q_0 x - p_0\right| = \left|q_0 \frac{p_1}{q_1} - p_0\right| = \frac{1}{q_1}\left|q_0 p_1 - q_1 p_0\right| \overset{(H)}{=} \frac{1}{q_1} \overset{(G)}{=} \frac{1}{a_1} < \frac{1}{2}$$

where (G) applies the recursion theorem and (H) the sign theorem. Because of (§§§) it is $\left|y_0 x - z_0\right| \geq \frac{1}{2}$, i.e. together $\left|q_0 x - p_0\right| < \left|y_0 x - z_0\right|$ which contradicts the definition of $y_0$ and $z_0$! This proves claim 1 for $N = 1$. □(claim 1)

Next we observe

<u>Claim 3</u>: $\frac{z_0}{y_0}$ is a best approximation of the second kind of x.

<u>Otherwise</u>: $\left|bx - a\right| \leq \left|y_0 x - z_0\right|$ for an $\frac{a}{b} \neq \frac{z_0}{y_0}$ with $b \leq y_0$, which contradicts the definition of $y_0$ and $z_0$! □(claim 3)

According to theorem 30, $\frac{z_0}{y_0}$ is a convergent of x, i.e. $\frac{z_0}{y_0} = \frac{p_s}{q_s}$ for an $s \leq k$. If $s = k$, the proof is done. Thus, we assume $s < k$.

<u>Claim 4</u>: For $s < k$, it is $\frac{1}{q_s + q_{s+1}} \geq \frac{1}{q_k + q_{k-1}}$

<u>Proof</u>: $s < k \Rightarrow s \leq k - 1 \Rightarrow q_s \leq q_{k-1}$ (corollary 5: denominators are monotonically increasing). Similarly, $s < k \Rightarrow s + 1 \leq k \Rightarrow q_{s+1} \leq q_k$. Together, this implies $q_k + q_{k-1} \geq q_s + q_{s+1}$. □(claim 4)

Next we get

$$\left|q_s x - p_s\right| = q_s \left|x - \frac{p_s}{q_s}\right| \overset{(I)}{>} q_s \frac{1}{\left(q_s + q_{s+1}\right) q_s} = \frac{1}{q_s + q_{s+1}} \overset{(J)}{\geq} \frac{1}{q_k + q_{k-1}}$$

where (I) is lemma 24 (lower bounds) and (J) is claim 4.

Furthermore, $\left|q_k x - p_k\right| = q_k \left|x - \frac{p_k}{q_k}\right| \overset{(K)}{<} q_k \frac{1}{q_k q_{k+1}} = \frac{1}{q_{k+1}}$, where (K) holds because of lemma 20 (upper bounds).

With $\frac{z_0}{y_0} = \frac{p_s}{q_s}$ and the definition of $y_0$ ($= q_s$) and $z_0$ ($= p_s$) (i.e. the minimizing property) it is $\left|q_s x - p_s\right| = \left|y_0 x - z_0\right| \leq \left|q_k x - p_k\right| \Rightarrow \frac{1}{q_k + q_{k-1}} \leq \frac{1}{q_{k+1}}$, which implies $q_{k+1} < q_k + q_{k-1}$: contradiction, because of the recursion theorem it is $q_{k+1} = a_{k+1} q_k + q_{k-1} \overset{(L)}{\geq} q_k + q_{k-1}$, where (L) holds with $a_k \geq 1$. Thus, $s = k$ which proves the overall theorem.

■

Putting the last two theorems together yields in

**Corollary 32**

$a/b$ is a best approximation of the second kind of x ⇔ x is a convergent of x.

■

According to theorem 31, every convergent is a best approximation of the second kind, and each best approximation of the second kind is also a best approximation of the first kind (lemma 28). We keep this observation as

**Note 33**
Every convergent is a best approximations of the first kind.

■

But are best approximations of the first kind also always convergents? Not quite: the next theorem proves that a best approximation of the first kind is a convergent or a semiconvergent.

**Theorem 34** (*Lagrange, 1798 —*
*1st Kind Best Approximations are Convergents or Semiconvergents*)
Let $a/b$ be a best approximation of the first kind of $x = \left[a_0; a_1, \cdots, a_N\right]$. Then is $a/b$ a convergent or a semiconvergent of x.

Proof: By definition it is $\left|x - \frac{c}{d}\right| > \left|x - \frac{a}{b}\right|$ for $\frac{c}{c} \neq \frac{a}{b}$ and $d \leq b$.

Claim 1: $a/b > a_0$

Otherwise: $\frac{a}{b} \leq a_0 = \frac{a_0}{1}$, thus $x - a_0 \leq x - \frac{a}{b}$. Now $x - a_0 = \frac{1}{a_1 + \ddots} > 0$ , thus $0 < x - a_0 \leq x - \frac{a}{b} \Rightarrow \left|x - \frac{a_0}{1}\right| \leq \left|x - \frac{a}{b}\right|$. Because $1 \leq b$ we got a contradiction since $a/b$ is a best approximation of the first kind. $\square_{\text{(claim 1)}}$

Claim 2: $a/b < a_0 + 1$

Otherwise: $\frac{a}{b} \geq a_0 + 1$ and based on the geometric situation depicted in figure 6, it follows that $\left|x - \frac{a_0 + 1}{1}\right| \leq \left|x - \frac{a}{b}\right|$ with $1 \leq b$, which contradicts $a/b$ being a best approximation of the first kind. $\square_{\text{(claim 2)}}$

Consequently, $a/b$ lies between $x_0 = a_0$ and $x_{-1,1} = a_0 + 1$ (see equation 22), i.e.

$$x_0 = a_0 < \frac{a}{b} < a_0 + 1 = x_{-1,1} \qquad (26)$$

and is, thus, covered by the set of intervals defined by the convergents and semiconvergents of x (see figure 6).

Assumption: $a/b$ is neither a convergent nor a semiconvergent.

This results in the following cases:

- Case 1: $a/b$ lies between two semiconvergents $x_{k-1,r}$ and $x_{k-1,r+1}$
- Case 2: $a/b$ lies between two convergents $x_k$ and $x_{k+2}$
- Case 3: $a/b$ lies between a convergent and a semiconvergent

We will show that all three cases lead to a contradiction, i.e. the assumption must be false, thus, the theorem is proven.

Case 1: $a/b$ lies between $x_{k-1,r} = \dfrac{rp_k + p_{k-1}}{rq_k + q_{k-1}}$ and $x_{k-1,r+1} = \dfrac{(r+1)p_k + p_{k-1}}{(r+1)q_k + q_{k-1}}$

Then,

$$\left|\frac{a}{b} - \frac{rp_k + p_{k-1}}{rq_k + q_{k-1}}\right| < \left|\frac{(r+1)p_k + p_{k-1}}{(r+1)q_k + q_{k-1}} - \frac{rp_k + p_{k-1}}{rq_k + q_{k-1}}\right| \overset{(A)}{=} \frac{1}{\left((r+1)q_k + q_{k-1}\right)\left(rq_k + q_{k-1}\right)}$$

where (A) results from the same computation performed in the proof of lemma 23. Furthermore, it is

$$\left|\frac{a}{b} - \frac{rp_k + p_{k-1}}{rq_k + q_{k-1}}\right| = \frac{\left|a\left(rq_k + q_{k-1}\right) - b\left(rp_k + p_{k-1}\right)\right|}{b\left(rq_k + q_{k-1}\right)} \overset{(B)}{\geq} \frac{1}{b\left(rq_k + q_{k-1}\right)} \quad (*)$$

where (B) is seen to be valid as follows: $a\left(rq_k + q_{k-1}\right) - b\left(rp_k + p_{k-1}\right) \in \mathbb{Z}$ and, thus, $\left|a\left(rq_k + q_{k-1}\right) - b\left(rp_k + p_{k-1}\right)\right| \in \mathbb{N}_0$; if it would be zero, the first modulus in (*) would be zero, i.e. $a/b = x_{k-1,r}$ which contradicts the assumption of the claim, which in turn implies $\left|a\left(rq_k + q_{k-1}\right) - b\left(rp_k + p_{k-1}\right)\right| \geq 1$.

Together,

$$\frac{1}{b\left(rq_k + q_{k-1}\right)} < \frac{1}{\left((r+1)q_k + q_{k-1}\right)\left(rq_k + q_{k-1}\right)} \Rightarrow \frac{1}{b} < \frac{1}{(r+1)q_k + q_{k-1}},$$

thus,

$$b > (r+1)q_k + q_{k-1} \quad (\S)$$

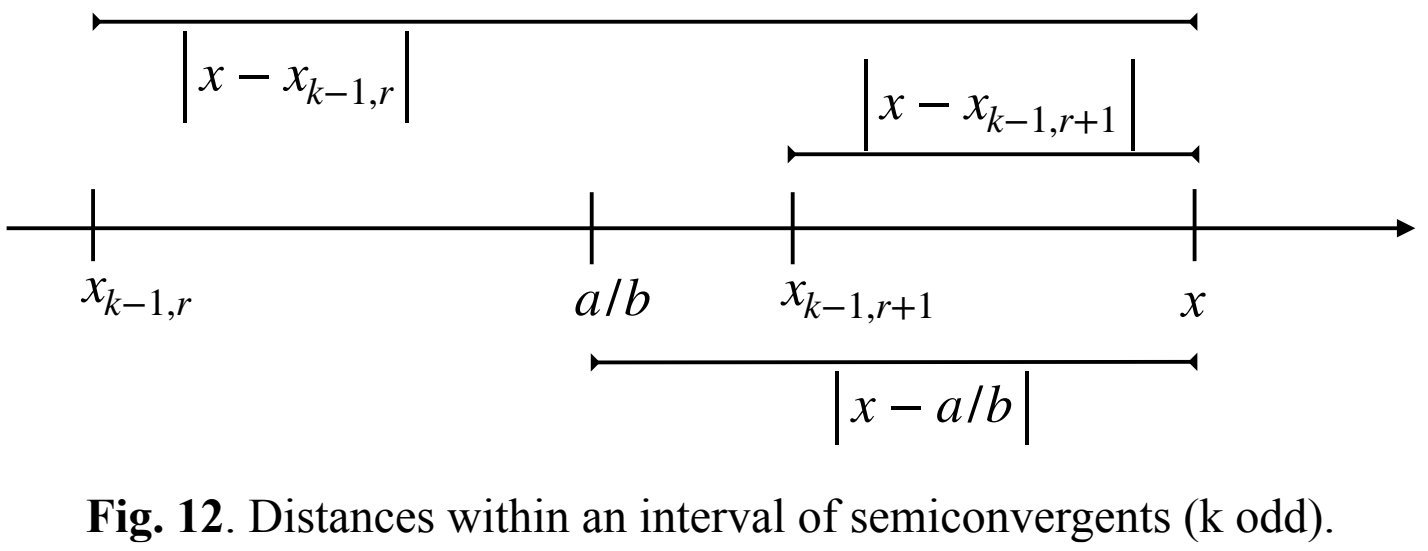

**Fig. 12**. Distances within an interval of semiconvergents (k odd).

Because of the monotony of the sequence of semiconvergents $(x_{s,t})_t$ (lemma 23), it is for an odd $k$ (i.e. $k-1$ even) $x_{k-1,r} < x_{k-1,r+1}$ (see the geometric situation in figure 12), thus

$$\left|x - \frac{a}{b}\right| > \left|x - \frac{(r+1)p_k + p_{k-1}}{(r+1)q_k + q_{k-1}}\right|$$

But with (§) it is $(r+1)q_k + q_{k-1} < b$, thus, $a/b$ is not a best approximation of the first kind to x, which is a contradiction. $k$ even leads to a contradiction too, i.e. case (1) is not possible $\square_{(\text{case } 1)}$

Case 2: $a/b$ lies between $x_k$ and $x_{k+2}$

Then, $\left|\frac{a}{b} - \frac{p_k}{q_k}\right| < \left|\frac{p_k}{q_k} - \frac{p_{k+2}}{q_{k+2}}\right| \overset{(C)}{=} \frac{a_{k+2}}{q_k q_{k+2}} < \frac{1}{q_k q_{k+2}}$ where (C) is equation (11) from corollary 12, and with note 29 it is $\left|\frac{a}{b} - \frac{p_k}{q_k}\right| \geq \frac{1}{b q_k}$.

Together, $\frac{1}{b q_k} < \frac{1}{q_k q_{k+2}} \Rightarrow \frac{1}{b} < \frac{1}{q_{k+2}} \Rightarrow b > q_{k+2}$. Because of the geometric situation shown in figure 13, it is $\left|x - \frac{a}{b}\right| > \left|x - \frac{p_{k+2}}{q_{k+2}}\right|$, which is a contradiction to $a/b$ being a best approximation of the first kind to x and $b > q_{k+2}$. $\square_{(\text{case } 2)}$

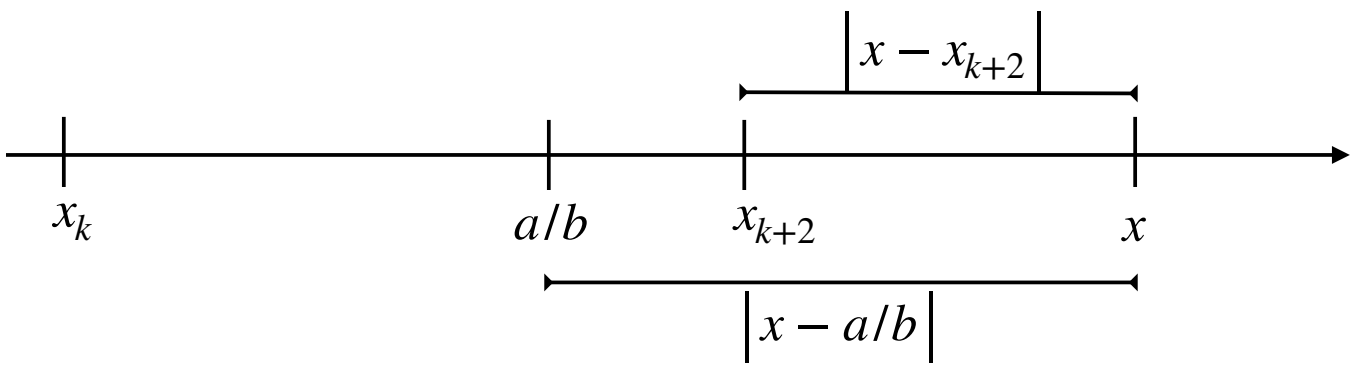

**Fig. 13**. Distances within an interval of convergents (k even).

Case 3: $a/b$ lies between a convergent and a semiconvergent.

This implies that $a/b$ lies between $x_k$ and $x_{k,1}$ (see figure 6), otherwise $a/b$ would lie between two semiconvergents that has already been covered in case 1.

Thus, $\left|\frac{a}{b} - \frac{p_k}{q_k}\right| < \left|x_k - x_{k,1}\right|$, but

$$\left|x_k - x_{k,1}\right| = \left|\frac{p_k}{q_k} - \frac{p_{k+1}+p_k}{q_{k+1}+q_k}\right| = \left|\frac{p_k(q_{k+1}+q_k) - q_k(p_{k+1}+p_k)}{q_k(q_{k+1}+q_k)}\right|$$

$$= \left|\frac{p_k q_{k+1} - q_k p_{k+1}}{q_k(q_{k+1}+q_k)}\right| \overset{(D)}{=} \frac{1}{q_k(q_{k+1}+q_k)}$$

where (D) is the sign theorem. I.e. it is $\left|\frac{a}{b} - \frac{p_k}{q_k}\right| < \frac{1}{q_k(q_{k+1}+q_k)}$. As before, with note 29 it is $\left|\frac{a}{b} - \frac{p_k}{q_k}\right| \geq \frac{1}{b q_k} \Rightarrow \frac{1}{b q_k} < \frac{1}{q_k(q_{k+1}+q_k)} \Rightarrow b > q_{k+1} + q_k$.

The geometric situation from figure 14 reveals $\left|x - \frac{a}{b}\right| > \left|x - \frac{p_k + p_{k-1}}{q_k + q_{k-1}}\right|$, which is a contradiction to $a/b$ being a best approximation of the first kind to x and $b > q_{k+1} + q_k$. $\square_{\text{(case 3)}}$

■

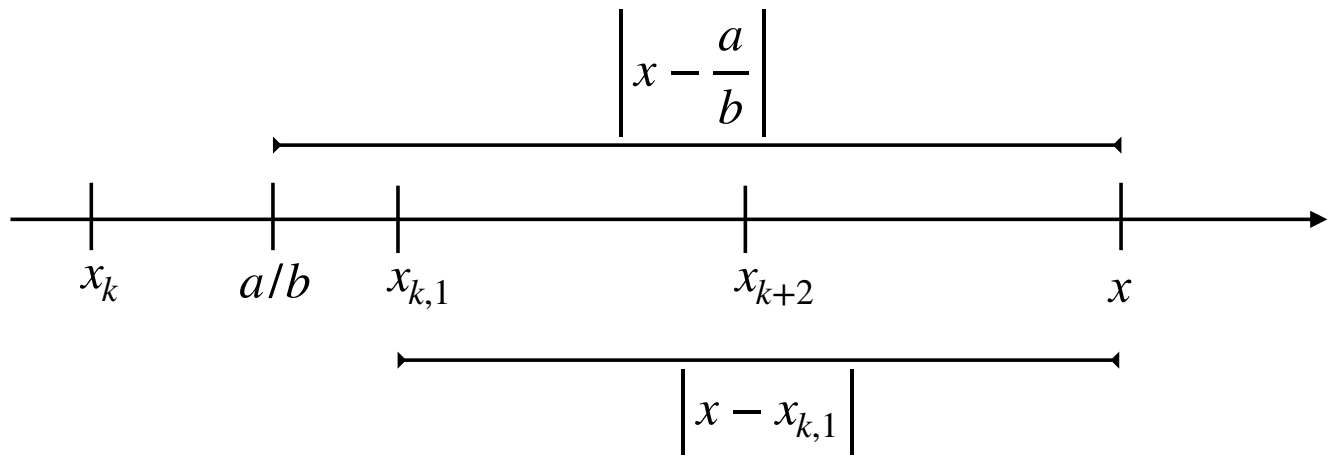

**Fig. 14**. Situation in which $a/b$ is between a convergent and its first semiconvergent (k even).

Finally, we give a simple criterion that allows to prove that a given fraction is a convergent of another real number. This theorem is a cornerstone of computing a prime factor with Shor's algorithm.

**Theorem 35** (*Legendre, 1798 — Convergent Criterion*)
Let $\left|x - \frac{a}{b}\right| < \frac{1}{2b^2} \Rightarrow a/b$ is a convergent of x.

Proof: We show that $a/b$ is a best approximation of the second kind of x. Theorem 30 then proves the claim.

Let $\left|dx - c\right| \leq \left|bx - a\right|$ for $\frac{a}{b} \neq \frac{c}{d}$ and $d > 0$. We need to prove $d > b$.

Now, $\left|bx - a\right| = b\left|x - \frac{a}{b}\right| < b\frac{1}{2b^2} = \frac{1}{2b}$. This implies $\left|dx - c\right| < \frac{1}{2b} \Leftrightarrow d\left|x - \frac{c}{d}\right| < \frac{1}{2b} \Leftrightarrow \left|x - \frac{c}{d}\right| < \frac{1}{2db}$. Thus,

$$\left|\frac{c}{d} - \frac{a}{b}\right| = \left|\frac{c}{d} - x + x - \frac{a}{b}\right| \leq \left|\frac{c}{d} - x\right| + \left|x - \frac{a}{b}\right| < \frac{1}{2db} + \frac{1}{2b^2} = \frac{b+d}{2db^2}$$

With note 29 (distance of fractions), it is also $\left|\frac{c}{d} - \frac{a}{b}\right| \geq \frac{1}{db}$. Together, it is

$$\frac{1}{db} < \frac{b+d}{2db^2} \Leftrightarrow 1 < \frac{b+d}{2b} \Leftrightarrow 2b < b + d \Leftrightarrow d > b.$$

■

# Part II: Probability of the Occurrence of Convergents

## 7. Estimating Secant Lengths

In this part, we use the main arguments of [5].

In order to estimate the probability of the occurrence of a certain state after having performed the quantum Fourier transform we need the following estimation of a lower bound and an upper bound of the length of a secant of the unit circle.

**Lemma 36** (*Secant Length Estimation*)

If $\varphi \in [-\pi, \pi]$ then $\frac{2\left|\varphi\right|}{\pi} \leq \left|1 - e^{i\varphi}\right| \leq \left|\varphi\right|$

Proof: The upper bound follows from elementary geometry, namely that the length of a secant is less than or equal to the length of the corresponding arc of a circle (see figure 15).

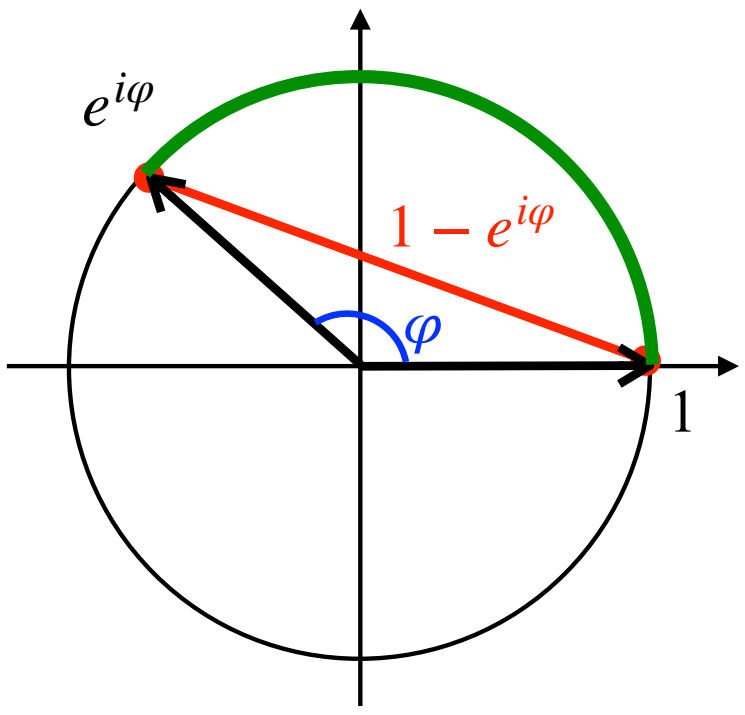

**Fig. 15**. The length of a secant is small than the arc of the corresponding unit circle.

The length of the arc determined by the angle $\varphi$ on a circle of radius r is $r\varphi$; i.e. if the circle is a unit circle, the length of the arc (green in the figure) is $\varphi$.

A secant of the unit circle (red in the figure) can be defined by the two complex numbers on the unit circle (black in the figure) that are the endpoints of the secant. Thus, the length of this secant is the difference of these complex numbers. One of these points can always be 1 because a corresponding rotation is length preserving; the other point is then $e^{i\varphi}$, where $\varphi$ is the angle of the arc cut by the secant. The length of this secant is then $\left|1 - e^{i\varphi}\right|$.

This proves the inequality $\left|1 - e^{i\varphi}\right| \leq \left|\varphi\right|$. □(upper bound)

Next, we compute

$$\left|1-e^{i\varphi}\right| \overset{(A)}{=} \left|1-\cos\varphi - i\sin\varphi\right| \overset{(B)}{=} \sqrt{\left(1-\cos\varphi\right)^2+\sin^2\varphi}$$
$$= \sqrt{1-2\cos\varphi+\cos^2\varphi+\sin^2\varphi} = \sqrt{2-2\cos\varphi}$$
$$= \sqrt{2}\sqrt{1-\cos\varphi} \overset{(C)}{=} \sqrt{2}\sqrt{2\sin^2\frac{\varphi}{2}}$$
$$\overset{(D)}{=} 2\sin\frac{\varphi}{2} \qquad (\S)$$

where (A) uses Euler's formula, (B) is the definition of the modulus of a complex number with $\mathrm{Re} = 1-\cos\varphi$ and $\mathrm{Im} = -\sin\varphi$, (C) is the double-angle formula, and (D) assumes that $\sin\frac{\varphi}{2} \geq 0$.

To estimate a lower bound for $\sin\frac{\varphi}{2}$, we analyze the function $f(x) = \sin x - \frac{2x}{\pi}$. From elementary calculus it is known that a function $\psi$ is concave on $D \subseteq \mathbb{R}$ if and only if its second derivative is not positive on D, i.e. $\psi'' \leq 0$ on D.

(Reminder: $\psi$ is *concave* on D :⇔

$$\forall x,y \in D\ \forall t \in [0,1] : \psi\left(tx+(1-t)y\right) \geq t\psi(x)+(1-t)\psi(y)$$

i.e. for any two points on the graph of $\psi$, the secant between these points is below the graph)

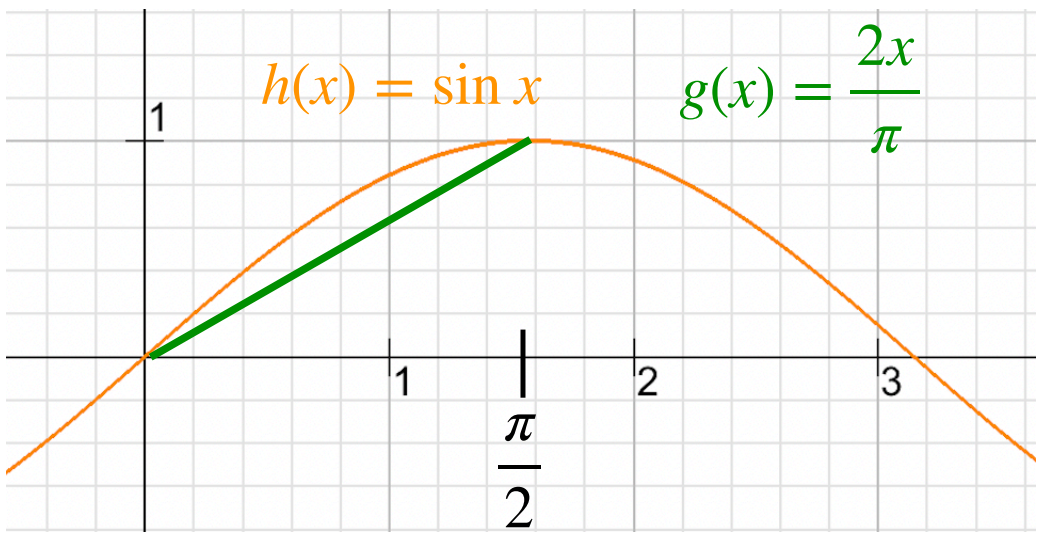

**Fig. 16**. The graphs of sin x and 2x/π .

With $d^2\sin x/dx^2 = -\sin x \leq 0$ especially for $x \in \left[0,\pi/2\right]$, i.e. $\sin x$ is concave on $x \in \left[0,\pi/2\right]$. Thus, the secant between $\sin 0 = 0$ and $\sin\frac{\pi}{2} = 1$ is below the graph of $\sin x$ (orange in figure 16). But this secant is given by $g(x) = \frac{2}{\pi}x$ (green in figure 16). Thus, it is $\frac{2}{\pi}x \leq \sin x$ for $x \in \left[0,\pi/2\right]$, i.e. with $\varphi := 2x$ we get $\frac{\varphi}{\pi} \leq \sin\frac{\varphi}{2}$ for $\varphi \in \left[0,\pi\right]$, and this implies $2\frac{\varphi}{\pi} \leq 2\sin\frac{\varphi}{2}$ for $\varphi \in \left[0,\pi\right]$.

Now, $\left|1-e^{i\varphi}\right| = 2\sin\frac{\varphi}{2}$ (see equation (§) above) implies $\left|1-e^{i\varphi}\right| \geq 2\frac{\varphi}{\pi}$ for $\varphi \in \left[0,\pi\right]$.

Furthermore, $\sin x$ is convex on $x \in \left[-\pi/2, 0\right]$ , thus an argument analogous to above shows that $\left|1 - e^{i\varphi}\right| \geq 2\frac{\left|\varphi\right|}{\pi}$ for $\varphi \in [-\pi, \pi]$ . $\square_{\text{(lower bound)}}$

■

## 8. Estimating Amplitude Parameters

As reminded in the introduction, the quantum part of Shor's algorithm produces in its final step the following quantum state via a measurement:

$$|y\rangle = \frac{1}{\sqrt{NA}} \sum_{j=0}^{A-1} \omega_N^{jpy} |y\rangle \tag{27}$$

Thus, according to the Born rule, the probability $P(y)$ of this particular state |y⟩ is the square of the modulus of the amplitude of |y⟩, i.e.

$$P(y) = \left| \frac{1}{\sqrt{NA}} \sum_{j=0}^{A-1} \omega_N^{jpy} \right|^2 = \frac{1}{NA} \left| \sum_{j=0}^{A-1} \omega_N^{jpy} \right|^2 \tag{28}$$

The argument of the modulus is a geometric sum $\sum q^j$ with $q = \omega_N^{py} = e^{\frac{2\pi i}{N} py}$, thus, in case $q \neq 1$,

$$P(y) = \frac{1}{NA} \left| \sum_{j=0}^{A-1} q^j \right|^2 = \frac{1}{NA} \left| \frac{1 - q^A}{1 - q} \right|^2 \tag{29}$$

with $q^A = e^{\frac{2\pi i}{N} Apy}$. In this section, in order to compute a lower bound for $P(y)$ we investigate some relations between the following parameters appearing in equation (28):

- n : the number to be factorized
- N : a power of 2 (e.g. $N = 2^m$) with $n^2 < N < 2n^2$
  - the choice of N effectively determines the domain of numbers that can be represented in the $|a\rangle$-part of the quantum register (see equation (1))
- p : the period of the modular exponentiation function $f(x) = a^x \bmod n$
- A : the number of arguments mapped to a given value of $f$

We also estimate bounds of the argument $2\pi \frac{Apy}{N}$ of $q^A = e^{\frac{2\pi i}{N} Apy}$.

### 8.1. Basics from Number Theory

For convenience, we remind the definition of the modulo function.

**Definition 37**:
The *modulo function* is the following map:

$$\begin{aligned} \text{mod} : \mathbb{N}_0 \times \mathbb{N} &\to \mathbb{N} \\ (z, n) &\mapsto z - \left\lfloor \frac{z}{n} \right\rfloor n \overset{\text{def}}{=} z \bmod n \end{aligned} \tag{30}$$

□

$z \bmod n$ is, thus, the residue left when dividing $z$ by $n$. I.e. if $r = z \bmod n$ then there is a number $k \in \mathbb{N}_0$ such that $z = kn + r$ with $0 \leq r < n$.

If $z \bmod n = \tilde{z} \bmod n = r$, we find numbers $k_1$ and $k_2$ such that $z = k_1 n + r$ and $\tilde{z} = k_2 n + r$ with $0 \leq r < n$. This implies that $z - \tilde{z} = (k_1 - k_2)n =: kn$, i.e. n is a divisor of $z - \tilde{z}$ (in symbols $n \,|\, (z - \tilde{z})$). We also get that $z \bmod n = \tilde{z} \bmod n$ implies that $\tilde{z} = z + kn$.

The equation $z \bmod n = \tilde{z} \bmod n$ is abbreviated as $z \equiv \tilde{z} \pmod{n}$ - in words: $z$ is *congruent* $\tilde{z}$ modulo n. As shown just before, $z \equiv \tilde{z} \pmod{n}$ is equivalent to $n \,|\, (z - \tilde{z})$ and to $\tilde{z} = z + kn$. We keep this as

**Note 38**
$z \equiv \tilde{z} \pmod{n} \Leftrightarrow n \,|\, (z - \tilde{z}) \Leftrightarrow \tilde{z} = z + kn$

■

Also, we remind the definition of modular exponentiation which turns out to play a key role in finding factors.

**Definition 39**:
For $0 < a < n$, the *modular exponentiation function* is the following map:

$$\begin{aligned} f : \mathbb{N}_0 &\to \mathbb{N}_0 \\ x &\mapsto a^x \bmod n \end{aligned} \tag{31}$$

□

The smallest number p that satisfies $f(x) = f(x + p)$ for all x is called the *period* of $f$. Especially, with $x = 0$ we get $f(0) = f(p)$ which means that $a^p \bmod n = a^0 \bmod n = 1 \bmod n$, i.e. $a^p \equiv 1 \pmod{n}$ which in turn is equivalent to $n \,|\, (a^p - 1)$. Thus, we have proven:

**Note 40**
$f(x) = a^x \bmod n$ has period p $\Leftrightarrow a^p \equiv 1 \pmod{n} \Leftrightarrow n \,|\, (a^p - 1)$

■

Finding a factor of n can be achieved by finding the period p of the function $f(x) = a^x \bmod n$. This is seen as follows: Let p be the period of $f$, then $n \,|\, (a^p - 1)$, i.e. $(a^p - 1) = kn$. Assume p is even (if p is odd, the algorithm of Shor is repeated with a different $a$, until an even p is found). With such an even p it is $\left(a^p - 1\right) = \left(a^{p/2} - 1\right)\left(a^{p/2} + 1\right) = kn$ which implies that $\left(a^{p/2} - 1\right)$ or $\left(a^{p/2} + 1\right)$ have a common divisor , which in turn means that gcd $\left(a^{p/2} - 1, n\right)$ or $\gcd\left(a^{p/2}+1, n\right)$ is a divisor of n. Thus, if an even period has been determined, classically efficient calculations can be used to compute a factor of n. If this factor is a prime number we can finish. Otherwise, we continue determining a factor of the former factor, and so on, until we end up with a prime factor of n.

Next, we determine an upper bound of the period p of the modular exponentiation by using group theory. A simple calculation shows that "≡" is an equivalence relation on $\mathbb{Z}$. The equivalence class of $z \in \mathbb{Z}$ is denoted as $[z]$ and is referred to as *residue class* of z modulo n. It is $[z] = \{\tilde{z} \in \mathbb{Z} \mid \tilde{z} \equiv z \pmod{n}\} = \{z + kn \mid k \in \mathbb{Z}\}$ (see note 38) where the latter set is sometimes written as $z + n\mathbb{Z}$. The set of all residue classes modulo n is denoted as $\mathbb{Z}_n$, i.e. $\mathbb{Z}_n = \left\{[0], [1], \ldots, [n-1]\right\}$.

We can multiply two residue classes modulo n as follows: $[x] \cdot [y] = [x \cdot y]$. With this multiplication, $\mathbb{Z}_n^* = \{[z] \in \mathbb{Z}_n \mid \gcd(z, n) = 1\}$ becomes a group. Because $\mathbb{Z}_n^* \subseteq \mathbb{Z}_n$, it is $\varphi(n) := \text{card}\mathbb{Z}_n^* \le \text{card}\mathbb{Z}_n = n$. Since every integer is a divisor of itself, it is $\gcd(n, n) = n \ne 1$ (for $n \ge 2$), i.e. the cardinality of numbers co-prime to n is less than n: $n \ge 2 \Rightarrow \varphi(n) < n$.

The well-known theorem of Lagrange from group theory states that for a group G with $\text{card}\, G = m < \infty$ and for each $x \in G$, it is $x^m = e$ (e is the unit element of G) — see lemma 3.2.5 in [8], for example. Thus, for $x \in \mathbb{Z}_n^*$ it is $x^{\varphi(n)} = 1$, i.e. $x^{\varphi(n)} \equiv 1 \pmod{n}$. Since the period p is the smallest number with $x^p \equiv 1 \pmod{n}$, it follows that $p \le \varphi(n)$ and, thus, $p < n$.

Now, the assumption of Shor's algorithm is that $0 < a < n$ and that $\gcd(a, n) = 1$, which ensures that $[a] \in \mathbb{Z}_n^*$, thus $[a]^p \equiv 1 \pmod{n}$ and $p < n$:

**Lemma 41**
Let p be the period of $f(x) = a^x \bmod n$. Then: $p < n$

■

## 8.2. Intervals of Consecutive Multiples of the Period

The relation between N and p is depicted in figure 17: multiples of N are always contained in closed intervals defined by consecutive multiples of p, i.e. it may happen that a multiple of N coincides with a multiple of p.

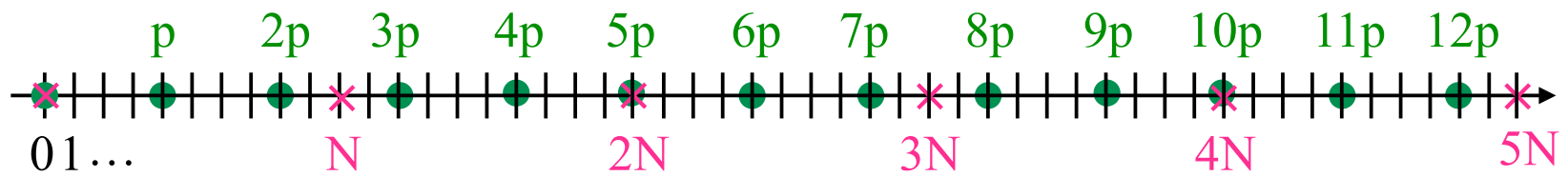

**Fig. 17**. Multiples of N are enclosed by immediately succeeding multiples of p.

**Note 42**
$\forall k \in \mathbb{N}\ \exists t \in \mathbb{N} : (t-1)p \le kN \le tp$

<u>Proof</u>: Pick an arbitrary $k \in \mathbb{N}$, i.e. $kN \in \mathbb{N}$ is also given. Then, we find a $\tilde{t} \in \mathbb{N}$ such that $\tilde{t}p \geq kN$ (trivial because $p > 0$). Let $t$ be the smallest of such $\tilde{t}$, i.e. $t \stackrel{\text{def}}{=} \min\left\{\tilde{t} \mid \tilde{t}p \geq kN\right\}$. Thus, $(t-1)p \leq kN$ because otherwise $(t-1)p > kN$ which is a contradiction because t was chosen minimal.

Together, $(t-1)p \leq kN \leq tp$. The claim follows because k an arbitrary number.

■

The situation we just discussed in shown in figure 18.

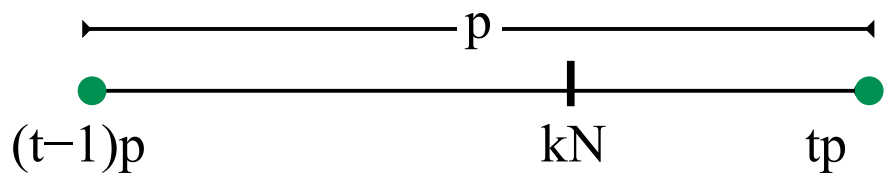

**Fig. 18**. Determining the interval of succeeding multiples of p enclosing a multiple of N.

Furthermore, two different multiples of N are in different intervals defined by succeeding multiples of p. Otherwise, the situation of figure 19 would occur implying that $N \leq p$ which is a contradiction as shown by the proof of note 43 below.

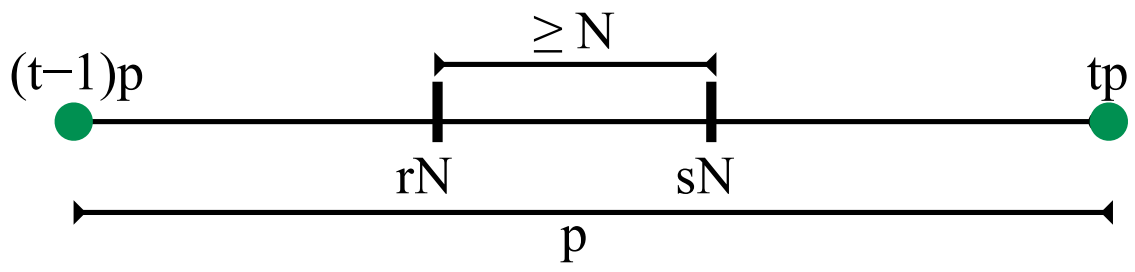

**Fig. 19**. No two multiples of N can be enclosed by succeeding multiples of p.

We denote a $t$ with $(t-1)p \leq kN \leq tp$ by $t_k$. This is justified by the next note 43 that proves that such a $t$ is uniquely determined by $k$. Especially, a multiple $kN$ is contained in “its” interval:

$$\forall k \in \mathbb{N}\ \exists t_k \in \mathbb{N} : kN \in \left[\left(t_k - 1\right)p, t_k p\right] \tag{32}$$

**Note 43**
Let $k \in \mathbb{N}$ and $t_k \in \mathbb{N}$ with $\left(t_k - 1\right)p \leq kN \leq t_k p$.

Then, $r \neq s \in \left\{0,\ldots,p-1\right\}$ implies $t_r \neq t_s$.

<u>Proof</u> (by contradiction): Assume $r \neq s$ but $t_r = t_s \stackrel{\text{def}}{=} t$ with $(t-1)p \leq rN \leq tp$ and $(t-1)p \leq sN \leq tp$ (see figure 19). W.l.o.g. $r < s \Rightarrow r+1 \leq s \Rightarrow sN - rN \geq (r+1)N - rN = N$. Also, $sN - rN \leq tp - (t-1)p = p$. Together, it is $N \leq sN - rN \leq p$, i.e. $N < p$.

According to lemma 41, we know $p < n \Rightarrow N < p < n < n^2$. But by selection of N (see the bullet list at the beginning of section 8) it is $n^2 < N$: contradiction.

■

The proof of note 43 has shown especially:

**Corollary 44**
$r \neq s \in \left\{0,\ldots,p-1\right\} \Rightarrow rN \notin \left[\left(t_s - 1\right)p, t_s p\right]$

■

By note 43, for $r \neq s \in \{0,...,p-1\}$ the numbers $t_r, t_s$ are different, i.e. for each $k \in \{0,...,p-1\}$ such a unique $t_k$ exists, i.e. the $p$ numbers $t_0, t_1, \ldots, t_{p-1}$ are different. Thus:

**Corollary 45**
There exist p different numbers $t_k$, $0 \leq k \leq p-1$, such that $(t_k - 1)p \leq kN \leq t_k p$.

■

These different numbers are strictly monotonically increasing.

**Note 46**
Let $k \in \mathbb{N}$ and $t_k \in \mathbb{N}$ with $(t_k - 1)p \leq kN \leq t_k p$. Then $t_k < t_{k+1}$.
Thus, $t_0 < t_1 < \ldots < t_{p-1}$.

Proof (by contradiction): Assume $t_{k+1} \leq t_k$, thus $t_{k+1} - 1 \leq t_k - 1$, which implies $t_{k+1}p \leq t_k p$ and $(t_{k+1} - 1)p \leq (t_k - 1)p$.
Now $kN < (k+1)N$, $(t_k - 1)p \leq kN \leq t_k p$ and $(t_{k+1} - 1)p \leq (k+1)N \leq t_{k+1}p$. This implies $(k+1)N \leq t_{k+1}p \leq t_k p$ and $(t_k - 1)p \leq kN < (k+1)N$, which finally results in $(t_k - 1)p < (k+1)N \leq t_k p$ — which is a contradiction to corollary 44 because this would imply that $(k+1)N \in \left[(t_k - 1)p, t_k p\right]$

■

Each multiple $kN$ of N is “close” to a multiple $tp$ in the sense that $kN$ is at most $p/2$ apart from $(t_k - 1)p$ or $t_k p$ (see figure 20).

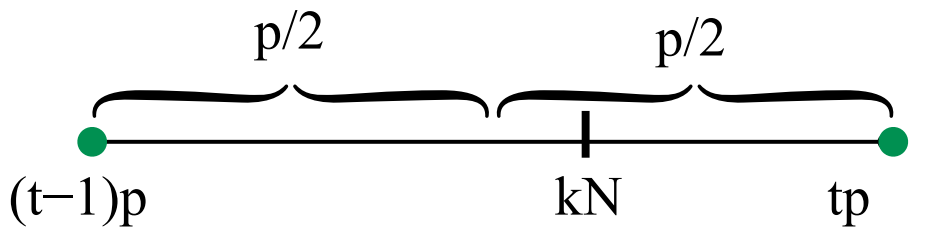

**Fig. 20**. A multiple of N is always “close” to a multiple of p.

More precisely:

**Note 47**
$\forall k \in \mathbb{N}\ \exists t \in \mathbb{N} : \left|(t-1)p - kN\right| \leq \frac{p}{2} \ \vee\ \left|tp - kN\right| \leq \frac{p}{2}$

Proof: It is $(t-1)p \leq kN \leq tp$, i.e. by definition $kN \in \left[(t-1)p, tp\right]$. This implies $kN - (t-1)p \leq \frac{p}{2} \vee tp - kN \leq \frac{p}{2}$ (see figure 20), otherwise:

$$kN - (t-1)p > \frac{p}{2} \wedge tp - kN > \frac{p}{2} \Leftrightarrow -(t-1)p > \frac{p}{2} - kN \wedge tp > \frac{p}{2} + kN$$

$\Rightarrow \; tp - (t-1)p > \frac{p}{2} - kN + \frac{p}{2} + kN = p$, but $tp - (t-1)p = p$, i.e. $p > p$: Contradiction! This proves the claim $\left|(t-1)p - kN\right| \leq \frac{p}{2} \;\vee\; \left|tp - kN\right| \leq \frac{p}{2}$.

■

And as before, from note 43 it follows that for $k \in \{0,..., p-1\}$ these numbers t are all different. Precisely:

**Corollary 48**
Let $0 \leq k \leq p-1$ and $t_k \in \mathbb{N}$ such that $\left|\left(t_k - 1\right)p - kN\right| \leq \frac{p}{2} \;\vee\; \left|t_k p - kN\right| \leq \frac{p}{2}$. If $r \neq s$, then $t_r \neq t_s$.

■

The multiples of N are sparsely scattered across the intervals of consecutive multiples of p. More precisely, intervals of consecutive multiples of p, which contain a multiple of N, are not consecutive. This is the content of

**Note 49**
Let $n > 2$, $k \in \{0,..., p-1\}$ and $t_k \in \mathbb{N}$ with $\left(t_k - 1\right)p \leq kN \leq t_k p$.

Then $t_{k+1} > t_k + 1$ as well as $t_{k-1} < t_k - 1$.

Proof: Because of note 46, it is $t_{k+1} > t_k$, thus, $t_{k+1} \geq t_k + 1$.

Assumption: $t_{k+1} = t_k + 1$
By definition: $(k+1)N \in \left[(t_{k+1} - 1)p, t_{k+1}p\right]$, and by assumption $t_{k+1} = t_k + 1$, thus it is $(k+1)N \in \left[t_k p, (t_k + 1)p\right]$. Also by definition, $kN \in \left[\left(t_k - 1\right)p, t_k p\right]$.

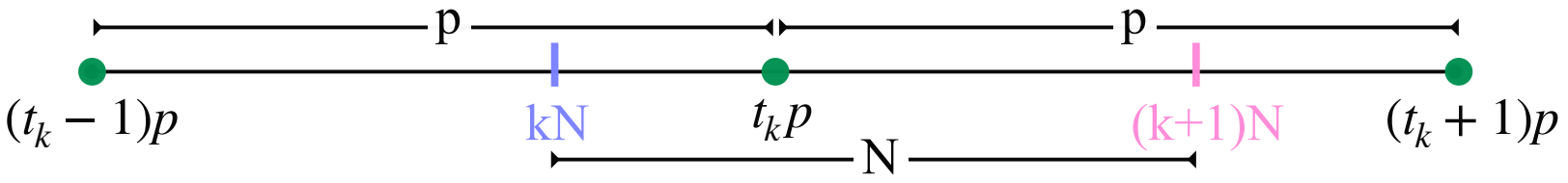

**Fig. 21**. Situation in case $kN$ and $(k+1)N$ lying within two consecutive intervals of consecutive multiples of p.

Now, $\left[\left(t_k - 1\right)p, t_k p\right], \left[t_k p, \left(t_k + 1\right)p\right] \subseteq \left[\left(t_k - 1\right)p, (t_k + 1)p\right]$ which implies $kN, (k+1)N \in \left[\left(t_k - 1\right)p, (t_k + 1)p\right]$.

Thus, $N = (k+1)N - kN \leq (t_k + 1)p - (t_k - 1)p = 2p$, i.e. $N \leq 2p$ (see figure 21).

By lemma 41, it is $p < n \Rightarrow 2p < 2n$. With $n > 2 \Rightarrow n^2 > 2n$ and by definition of N it is $n^2 < N$, thus $N > n^2 > 2n > 2p$: contradiction!

Thus, the assumption is false, which implies $t_{k+1} > t_k + 1$. The claim $t_{k-1} < t_k - 1$ is proven similarly.

■

The resulting geometric situation is depicted in figure 22.

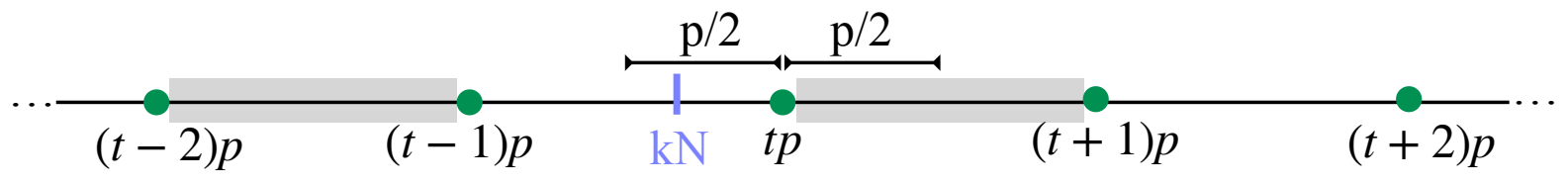

**Fig. 22**. If $kN$ is in an interval defined by two consecutive multiples of p, the preceding and succeeding intervals do not contain a multiple of N.

If $kN \in \left[(t_k - 1)p, t_k p\right]$ it is $t_k p - kN \le p/2$ or $kN - (t_k - 1)p \le p/2$ (see figure 20 or figure 22). In case $kN - (t_k - 1)p \le p/2$ we define $\hat{t} := t_k - 1$ and $kN - \hat{t}p \le p/2$ results, and in case of $t_k p - kN \le p/2$ we define $\hat{t} := t_k$ implying $\hat{t}p - kN \le p/2$. And according to note 43, this $\hat{t}$ is uniquely defined. This proves

**Note 50**
$\forall\, k \in \mathbb{N}\ \exists!\, \hat{t} \in \mathbb{N} : \left|\hat{t}p - kN\right| \le \frac{p}{2}$, i.e. $\hat{t}$ is uniquely determined by $k$.

■

This is next rewritten into a format more useful for what follows.

**Note 51**
Let $k \in \{0,\dots,p-1\}$ and $t_k \in \mathbb{N}$ with $t_k \in \left[k\frac{N}{p} - \frac{1}{2}, k\frac{N}{p} + \frac{1}{2}\right]$.
If $r \ne s \in \{0,\dots,p-1\}$ then $t_r \ne t_s$.

Proof: It is $t_k \in \left[k\frac{N}{p} - \frac{1}{2}, k\frac{N}{p} + \frac{1}{2}\right] \Leftrightarrow k\frac{N}{p} - \frac{1}{2} \le t_k \le k\frac{N}{p} + \frac{1}{2}$

$\Leftrightarrow kN - \frac{p}{2} \le pt_k \le kN + \frac{p}{2} \Leftrightarrow -\frac{p}{2} \le pt_k - kN \le \frac{p}{2} \Leftrightarrow \left|pt_k - kN\right| \le \frac{p}{2}$.
Note 50 shows that $t_k$ is uniquely determined by $k$.

■

Finally, we can prove the following

**Corollary 52**
There exist p different numbers $t_k$, $0 \le k \le p - 1$, such that
$t_k \in \left[k\frac{N}{p} - \frac{1}{2}, k\frac{N}{p} + \frac{1}{2}\right]$.

Proof: There exist p different numbers $t_k$, $0 \le k \le p - 1$, such that $(t_k - 1)p \le kN \le t_k p$ (corollary 45). The proof of note 47 shows that this implies $\left|(t_k - 1)p - kN\right| \le \frac{p}{2} \vee \left|t_k p - kN\right| \le \frac{p}{2}$. And the proof of note 51 shows that this implies $t_k \in \left[k\frac{N}{p} - \frac{1}{2}, k\frac{N}{p} + \frac{1}{2}\right]$.

■

### 8.3. Cardinality of Pre-Images

First we show that the parameter A is greater than 1, i.e. at least two numbers available in the $|a\rangle$-part of the quantum register are mapped by $f$ to the same value.

**Note 53**

$A > 1$

Proof: As reminded in the introduction, the quantum Fourier transform of the Shor algorithm produces the following state:

$$|a\rangle|b\rangle = \frac{1}{\sqrt{N}} \sum_{x=0}^{N-1} |x\rangle|f(x)\rangle$$

After measurement of the $|b\rangle$-part of the register, the $|a\rangle$-part is in the state

$$|a\rangle = \frac{1}{\sqrt{A}} \left( |x\rangle + |x+p\rangle + |x+2p\rangle + \cdots + |x+(A-1)p\rangle \right) \quad (33)$$

i.e. $f^{-1}(|x\rangle) = \left\{ |x\rangle, |x+p\rangle, |x+2p\rangle, \cdots, |x+(A-1)p\rangle \right\}$.

Choose $x < p$ — such an x exists because otherwise it would be $p = 0$ but a period p satisfies $p > 0$. With of $p < n$ (lemma 41) and $n < n^2 < N$ (by choice of N) it is $x + p < 2p < N$ (see the proof of note 49). Thus, $x + p$ is in the domain of $f$ (in the sense that it is a value in the $|a\rangle$-part of the quantum register available as an argument for $f$), i.e. $f(x + p)$ is available in the $|b\rangle$-part of the register.

$A = 1$ would imply that $f^{-1}(|x\rangle) = \left\{ |x\rangle \right\}$ and, thus, $|x+p\rangle \notin f^{-1}(|x\rangle)$, i.e. $f(|x\rangle) \neq f(|x+p\rangle)$ for $p \neq 0$. Since $p > 0$, the function $f$ would not be periodic.

■

Next, we prove tighter bounds for the parameter A.

**Note 54**

$(A-1)p < N < (A+1)p$

Proof: As in the proof of note 53, we choose $x < p$. With

$$|a\rangle = \frac{1}{\sqrt{A}} \left( |x\rangle + |x+p\rangle + |x+2p\rangle + \cdots + |x+(A-1)p\rangle \right),$$

i.e. $\left\{ |x\rangle, |x+p\rangle, |x+2p\rangle, \cdots, |x+(A-1)p\rangle \right\}$ are all values in the $|a\rangle$-part of the register being mapped to $f(x)$, i.e. $A$ is the largest number satisfying $x + (A-1)p < N$. With $x \geq 0$, this implies $(A-1)p < N$ — which is the first part of the claim.

Thus, $(A+1)p > N$. Otherwise $(A+1)p \leq N$ and with $x < p$ it would be $x + Ap < p + Ap = (A+1)p \leq N$, i.e $|x + Ap\rangle$ would also be in the $|a\rangle$-part of

the register being mapped to $f(x)$, which is a contradiction to the definition of A. This proves the second part of the claim.

■

The next estimation gives an approximation of $N$ in terms of the product of $A$ and $p$.

**Note 55**

$N \approx Ap$

Proof: Because of $(A-1)p < N < (A+1)p$, the geometric situation is as depicted in figure 23, i.e. $N \in [(A-1)p, Ap]$ or $N \in [Ap, (A-1)p]$. Thus $|N - Ap| \leq p$.

Now, $p < n$ (lemma 41) and $n < n^2 < N$ by choice of N. In practice, $n$ is a large number, i.e. $n^2$ is huge compared to $n : n \ll n^2 < N$. Together:

$$p \ll N \tag{34}$$

In this sense, $p$ is a small number, i.e. $|N - Ap|$ is small too: $N \approx Ap$.

■

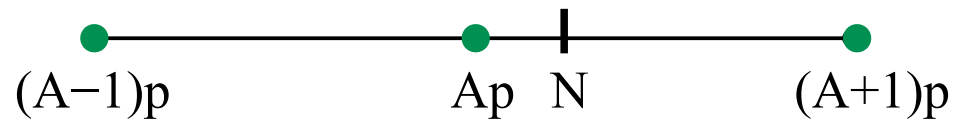

**Fig. 23**. $N$ is embraced by $(A-1)p$ and $(A+1)p$.

## 8.4. Estimating Arguments of Amplitudes of Potential Measurement Results

The next lemma is the main result of this section for what follows.

**Lemma 56**

Let $y \in \left[k\frac{N}{p} - \frac{1}{2}, k\frac{N}{p} + \frac{1}{2}\right]$ and $k \in \{0, ..., p-1\}$. Then:

$$2\pi\frac{y(A-1)p}{N}, 2\pi\frac{py}{N} \in [-\pi, +\pi]$$

Proof: It is $y \in \left[k\frac{N}{p} - \frac{1}{2}, k\frac{N}{p} + \frac{1}{2}\right] \Leftrightarrow k\frac{N}{p} - \frac{1}{2} \leq y \leq k\frac{N}{p} + \frac{1}{2} \Leftrightarrow$ (multiply with p) $kN - \frac{p}{2} \leq py \leq kN + \frac{p}{2} \Leftrightarrow -\frac{p}{2} \leq py - kN \leq \frac{p}{2} \Leftrightarrow$

$$yp - kN \leq \frac{p}{2} \;\wedge\; kN - yp \leq \frac{p}{2} \tag{§}$$

By note 54, it is $(A-1)p < N \Rightarrow \frac{(A-1)p}{N} < 1$ (§§).

Furthermore: $\frac{(A-1)p}{N} = \frac{Ap}{N} - \frac{p}{N} \overset{(A)}{\approx} \frac{N}{N} - \frac{p}{N} = 1 - \frac{p}{N} \overset{(B)}{\approx} 1$ (§§§) ,

where (A) is because of note 55 ($N \approx Ap$), and (B) is because of equation (34) ($p \ll N$).

Next we compute the lower bound for the first fraction of the claim:

$$\begin{aligned} 2\pi \frac{y(A-1)p}{N} &= 2\pi \frac{(A-1)}{N} yp \\ &\overset{(C)}{\geq} 2\pi \frac{(A-1)}{N}\left(kN - \frac{p}{2}\right) \\ &= 2\pi(A-1)k - \pi \frac{(A-1)p}{N} \\ &\overset{(D)}{\geq} -\pi \frac{(A-1)p}{N} \overset{(E)}{\approx} -\pi \end{aligned}$$

where (C) follows from the second inequation of (§) above, (D) is because of $2\pi(A-1)k \geq 0$, and (E) is implied by (§§§) above.

The upper bound for the first fraction of the claim is computed next:

$$\begin{aligned} 2\pi \frac{y(A-1)p}{N} &= 2\pi \frac{(A-1)}{N} yp \overset{(F)}{\leq} 2\pi \frac{(A-1)}{N}\left(kN + \frac{p}{2}\right) \\ &= 2\pi(A-1)k + \pi \frac{(A-1)p}{N} \overset{(G)}{<} 2\pi(A-1)k + \pi \\ &\overset{(H)}{<} 2\pi(A-1)p + \pi \overset{(I)}{<} 2\pi N + \pi \overset{(J)}{\equiv} \pi \end{aligned}$$

where (F) is implied by the first inequation of (§) above, (G) is (§§) above, (H) follows from the prerequisite $k \in \{0,...,p-1\}$, i.e. $k < p$, and (I) is the first inequation of (§§) above. Finally, we will estimate $e^{i\varphi}$ and because of $e^{i2\pi N} = 1$, (J) is justified.

Together, $-\pi \leq 2\pi \frac{y(A-1)p}{N} \leq \pi$, which proves the first claim. $\square_{\text{(first fraction)}}$

Next,

$$\begin{aligned} 2\pi \frac{py}{N} &= \frac{2\pi}{N} py \overset{(K)}{\leq} \frac{2\pi}{N}\left(kN + \frac{p}{2}\right) \\ &= 2\pi k + \pi \frac{p}{N} \overset{(L)}{<} 2\pi k + \pi \overset{(M)}{\equiv} \pi \end{aligned}$$

with (K) from the first inequation of (§) before, (L) because $p < N$, and (M) because we will estimate $e^{i\varphi}$. I.e. the upper bound of the second fraction is as claimed.

The correctness of the lower bound is seen as follows:

$$\begin{aligned} 2\pi \frac{py}{N} &= \frac{2\pi}{N} py \overset{(N)}{\geq} \frac{2\pi}{N}\left(kN - \frac{p}{2}\right) \\ &= 2\pi k - \pi \frac{p}{N} \overset{(O)}{>} -\pi \frac{p}{N} \overset{(Q)}{>} -\pi \end{aligned}$$

with the second inequation of (§) before gives (N), (O) is because of $2\pi k > 0$, and (Q) is true because $0 < p < N$, thus $0 < p/N < 1$. $\square_{\text{(second fraction)}}$

■

## 9. Estimating Probabilities

We are now ready to compute the probability $P(y)$ that the state |y⟩, which is prepared by the quantum part of Shor's algorithm, is "close" (i.e. within a distance of $1/2$) to a multiple of $p/N$.

**Lemma 57**

Assume $q = e^{i2\pi\frac{yp}{N}} \neq 1$ and let P be the probability that $y \in \left[k\frac{N}{p} - \frac{1}{2}, k\frac{N}{p} + \frac{1}{2}\right]$ for a $k \in \{0,\ldots,p-1\}$. Then $P \approx \frac{4}{\pi^2}$.

<u>Proof</u>: According to equation (29) the probability $P(y)$ to measure a particular $y \in \left[k\frac{N}{p} - \frac{1}{2}, k\frac{N}{p} + \frac{1}{2}\right]$ is

$$P(y) = \frac{1}{NA}\left|\frac{1-q^A}{1-q}\right|^2 \tag{35}$$

in case $q \neq 1$ (which is the assumption) where $q = e^{i2\pi\frac{yp}{N}}$ (the case $q = 1$ will be treated separately in note 58). Thus, with $q^A = e^{i2\pi\frac{yAp}{N}}$ it is

$$\left|\frac{1-q^A}{1-q}\right| = \frac{\left|1-q^A\right|}{\left|1-q\right|} = \frac{\left|1-e^{i2\pi\frac{yAp}{N}}\right|}{\left|1-e^{i2\pi\frac{yp}{N}}\right|} \tag{36}$$

The structure of the numerator and denominator recommends to estimate both by means of the lemma 36 (secant length estimation). But applying lemma 36 requires that $2\pi\frac{yAp}{N}, 2\pi\frac{yp}{N} \in [-\pi, +\pi]$. By lemma 56 we know that under the prerequisite $y \in \left[k\frac{N}{p} - \frac{1}{2}, k\frac{N}{p} + \frac{1}{2}\right]$ it is $2\pi\frac{py}{N} \in [-\pi, +\pi]$, but unfortunately, lemma 56 only guarantees that $2\pi\frac{y(A-1)p}{N} \in [-\pi, +\pi]$.

Now, consider the following calculation:

$$\left|\frac{1-q^A}{1-q}\right| = \left|\frac{1-q^{A-1}}{1-q} + q^{A-1}\right| \overset{(A)}{\geq} \left|\frac{1-q^{A-1}}{1-q}\right| - \left|q^{A-1}\right| \overset{(B)}{=} \left|\frac{1-q^{A-1}}{1-q}\right| - 1 \quad (37)$$

where (A) holds because of $|a+b| \geq |a| - |b|$, and $\left|e^{i\varphi}\right| = 1$ implies (B): $\left|q^t\right| = \left|\left(e^{i2\pi\frac{yp}{N}}\right)^t\right| = \left|e^{i\left(2\pi\frac{ypt}{N}\right)}\right| = 1.$

Equation (37) allows to apply the secant length estimation (lemma 36) because in

$$\left|\frac{1-q^{A-1}}{1-q}\right| = \frac{\left|1-q^{A-1}\right|}{\left|1-q\right|} = \frac{\left|1-e^{i2\pi\frac{y(A-1)p}{N}}\right|}{\left|1-e^{i2\pi\frac{yp}{N}}\right|} \quad (38)$$

it is now $2\pi\frac{y(A-1)p}{N}, 2\pi\frac{yp}{N} \in [-\pi, +\pi]$ according to lemma 56.

First, we use the second inequation of $\frac{2\left|\varphi\right|}{\pi} \leq \left|1-e^{i\varphi}\right| \leq \left|\varphi\right|$ from lemma 36 with $\varphi = 2\pi\frac{yp}{N} \in [-\pi, +\pi]$ and get

$$\left|1-e^{i2\pi\frac{yp}{N}}\right| \leq 2\pi\frac{yp}{N} \quad (39)$$

Then, we use the first inequation of $\frac{2\left|\varphi\right|}{\pi} \leq \left|1-e^{i\varphi}\right| \leq \left|\varphi\right|$ from lemma 36 with $\varphi = 2\pi\frac{y(A-1)p}{N} \in [-\pi, +\pi]$ and get

$$\left|1-e^{i2\pi\frac{y(A-1)p}{N}}\right| \geq \frac{2}{\pi}\cdot 2\pi\frac{y(A-1)p}{N} \quad (40)$$

Using equations (39) and (40) in equation (38) results in

$$\left|\frac{1-q^{A-1}}{1-q}\right| = \frac{\left|1-e^{i2\pi\frac{y(A-1)p}{N}}\right|}{\left|1-e^{i2\pi\frac{yp}{N}}\right|} \geq \frac{2}{\pi}\cdot 2\pi y\frac{(A-1)p}{N}\cdot\frac{N}{2\pi yp} = \frac{2(A-1)}{\pi} \quad (41)$$

This result is now used in equation (37) (step (C) below) and we get

$$\left|\frac{1-q^A}{1-q}\right| = \left|\frac{1-q^{A-1}}{1-q}\right| - 1 \overset{(C)}{\geq} \frac{2(A-1)}{\pi} - 1$$
$$= \frac{2A}{\pi} - \frac{2}{\pi} - 1 = \frac{2A}{\pi} - \left(\frac{2}{\pi} + 1\right) \tag{42}$$

Using equation (42) in equation (35) (step (D) below) results in

$$P(y) = \frac{1}{NA}\left|\frac{1-q^A}{1-q}\right|^2 \overset{(D)}{\geq} \frac{1}{NA}\left(\frac{2A}{\pi} - \left(\frac{2}{\pi}+1\right)\right)^2$$
$$= \frac{1}{NA}\left(\frac{4A^2}{\pi^2} - \frac{4A}{\pi}\left(\frac{2}{\pi}+1\right) + \left(\frac{2}{\pi}+1\right)^2\right)$$
$$= \frac{1}{NA}\left(\frac{4A^2}{\pi^2} - \frac{8A}{\pi^2} - \frac{4A}{\pi} + \frac{4}{\pi^2} + \frac{4}{\pi} + 1\right)$$
$$= \frac{4A}{\pi^2 N} - \frac{8}{\pi^2 N} - \frac{4}{\pi N} + \frac{4}{\pi^2 NA} + \frac{4}{\pi NA} + \frac{1}{NA}$$
$$\geq \frac{4A}{\pi^2 N} - \frac{8}{\pi^2 N} - \frac{4}{\pi N} = \frac{4A}{\pi^2 N} - \frac{4}{\pi N}\left(1 + \frac{2}{\pi}\right)$$

Thus,

$$P(y) \geq \frac{4A}{\pi^2 N} - \frac{4}{\pi N}\left(1 + \frac{2}{\pi}\right) \tag{43}$$

According to note 55, we know $N \approx Ap \Rightarrow \frac{A}{N} \approx \frac{1}{p}$, i.e.

$$\frac{4A}{\pi^2 N} \approx \frac{4}{\pi^2 p} \tag{44}$$

Furthermore, since N is a “huge” number, we know that the following is “small”:

$$\frac{4}{\pi N}\left(1 + \frac{2}{\pi}\right) \overset{\text{def}}{=} \varepsilon \tag{45}$$

By using equations (44) and (45) in equation (43) results in

$$P(y) \geq \frac{4}{\pi^2}\frac{1}{p} - \varepsilon \tag{46}$$

for each $y \in \left[k\frac{N}{p} - \frac{1}{2}, k\frac{N}{p} + \frac{1}{2}\right]$. According to corollary 52, there exist p different numbers $y_k$ with $y_k \in \left[k\frac{N}{p} - \frac{1}{2}, k\frac{N}{p} + \frac{1}{2}\right]$ and for each of them is $P(y_k) \geq \frac{4}{\pi^2}\frac{1}{p} - \varepsilon$. Since we are not interested in a particular $y_k$ but in any of them, we need to sum up all probabilities $P(y_k)$ to get the overall probability $P$:

$$P = \sum_{i=0}^{p-1} P\left(y_i\right) \geq \frac{4}{\pi^2} - p\varepsilon \approx \frac{4}{\pi^2}$$

This proves the claim.

■

We still need to estimate the probability for the case $q = 1$.

**Note 58**

Let $q = 1$. Then $P(y) = \frac{A}{N}$.

Proof: In case $q = 1$,the probability is

$$P(y) = \frac{1}{NA}\left|\sum_{j=0}^{A-1} q^A\right|^2 = \frac{1}{NA}\left|\sum_{j=0}^{A-1} 1\right|^2 = \frac{1}{NA}A^2 = \frac{A}{N}.$$

■

## 10. Computing the Period

Let y be the result of the measurement produced by the Shor algorithm. Under the assumption that $q \neq 1$, the following holds:

**Theorem 59**

With probability $P \approx \frac{4}{\pi^2}$ there exists a $k \in \{0,...,p-1\}$, such that

$$\left|\frac{y}{N} - \frac{k}{p}\right| < \frac{1}{2p^2}$$

Proof: According to lemma 57, the probability that $y \in \left[k\frac{N}{p} - \frac{1}{2}, k\frac{N}{p} + \frac{1}{2}\right]$ for a $k \in \{0,...,p-1\}$ is $\approx 4/\pi^2$.

But $y \in \left[k\frac{N}{p} - \frac{1}{2}, k\frac{N}{p} + \frac{1}{2}\right] \Leftrightarrow -\frac{1}{2} \leq y - \frac{kN}{p} \leq +\frac{1}{2}$. Dividing the latter inequations by $N$ yields: $y \in \left[k\frac{N}{p} - \frac{1}{2}, k\frac{N}{p} + \frac{1}{2}\right] \Leftrightarrow -\frac{1}{2N} \leq \frac{y}{N} - \frac{k}{p} \leq +\frac{1}{2N}$.

Thus, $\left|\frac{y}{N}-\frac{k}{p}\right| \le \frac{1}{2N}$. By choice of $N$ it is $n^2 < N$. Furthermore, $p < n \Rightarrow p^2 < n^2 \Rightarrow p^2 < N \Rightarrow \frac{1}{N} < \frac{1}{p^2}$. This results in $\left|\frac{y}{N}-\frac{k}{p}\right| \le \frac{1}{2N} < \frac{1}{2p^2}$.

■

The theorem of Legendre (theorem 35) proves immediately:

**Theorem 60**
With probability $\approx 4/\pi^2$, $k/p$ is a convergent of $y/N$.

■

### 10.1. Determining the Period by Convergents: $q \neq 1$

The following algorithm determines with probability of approximately $4/\pi^2$ the period p we are looking for; is is applicable in the case $q \neq 1$:

1. Compute $\frac{y}{N} \in \mathbb{Q}_{>0}$
   The result of the measurement is $y \in \mathbb{N}$ and $N \in \mathbb{N}$ has been chosen
   $\Rightarrow \frac{y}{N} \in \mathbb{Q}_{>0}$ can be computed
2. Compute the continued fraction representation $[a_0; a_1, \ldots, a_m]$ of $\frac{y}{N} \in \mathbb{Q}$
3. Compute the convergents $[a_0; a_1, \ldots, a_u] = \frac{g_u}{h_u}$, $1 \le u \le m$
4. Determine $h_\omega$ with $h_\omega \ge h_u$ for $1 \le u \le m$ and $h_\omega < n$
   $\Rightarrow \frac{g_\omega}{h_\omega}$ is a very good approximation of $\frac{k}{p}$ because $\frac{1}{2h_\omega^2} \le \frac{1}{2h_u^2}$
5. Thus, $h_\omega \approx p$ is a candidate for the period p
6. Check whether p is in fact the period

### 10.2. Determining the Period by Convergents: $q = 1$

In case $q = 1$, the above algorithm is not applicable. But $q = 1 \Leftrightarrow e^{\frac{2\pi i}{N}py} = 1 \Leftrightarrow \frac{py}{N} \in \mathbb{Z} \Leftrightarrow p = k\frac{N}{y}$ with $k \in \mathbb{Z}$. Thus, the following algorithm can be used:

1. Compute $\frac{N}{y} \in \mathbb{Q}_{>0}$
   The result of the measurement is $y \in \mathbb{N}$ and $N \in \mathbb{N}$ has been chosen
   $\Rightarrow \frac{N}{y} \in \mathbb{Q}_{>0}$ can be computed

2. Select $k \in \mathbb{N}$
3. Compute $k \dfrac{N}{y}$
4. If $k \dfrac{N}{y} \notin \mathbb{N}$ go back to step (2)
5. If $k \dfrac{N}{y} \geq n$ go back to step (2)
6. $p = k \dfrac{N}{y}$ is a candidate for the period p
7. Check whether p is in fact the period
8. If p is not the period
   a. If some predefined termination criterion is met: stop
   b. Go back to step (2)

This may yield the period p but does not guarantee it.

## 11. Conclusion and Related Work

The literature analyzing, discussing, and refining Shor's algorithm [7] is vast. Of course, most text books on quantum computing explain the algorithm too (e.g. [4], [6]). In doing so, all this literature puts a sharp focus on the quantum part of the algorithm and sketches its classical parts at various depth. However, the mathematical treatment of the classical aspects is sketchy omitting most of the details and leave them as an exercise for the reader with references to corresponding text books from mathematics like [2] or [3]. The lecture notes by Preskill [5] go a bit deeper especially on the estimation of probabilities, but still omit the low-level details; however, the authors of the contribution at hand benefited a lot by the treatment in [5]. Finally, the anonymous text [1] refers to the arguments of [5] w.r.t. probability estimation and works them out. The authors of the contribution at hand build on this and provide many further details.

In doing so, the contribution at hand is very detailed on the probability estimation of being able to use the Legendre theorem in Shor's algorithm. The authors are not aware of any other publication providing these low-level details.

Furthermore, the contribution at hand is a self-contained treatment on continued fractions up to the Legendre theorem. All background is presented that is needed to understand this theorem, including all proofs with low-level details step-by-step.

The authors hope to foster the comprehension of the classical aspects of Shor's algorithm even at the level of beginners in quantum computing.

**Acknowledgements**

This work was partially funded by the BMWK project PlanQK (01MK20005N).